\documentclass[10pt]{article}

\usepackage{latexsym}
\usepackage{amsfonts}
\usepackage{amsmath}
\usepackage{mathrsfs}  
\usepackage{dsfont}    
\usepackage{bbold}     
\usepackage[english]{babel}
\usepackage{caption}
\usepackage{epsfig}
\usepackage{float}
\usepackage{subfigure}

\newtheorem{lemma}{Lemma}

\newtheorem{rmk}{Remark}

\newcommand{\zerarcounters}{\setcounter{equation}{0}}

\newcommand{\ZZZ}{\mathds{Z}}

\newcommand{\NNN}{\mathds{N}}
\newcommand{\QQQ}{\mathds{Q}}
\newcommand{\RRR}{\mathds{R}}

\newcommand{\one}{\mathds{1}}

\newcommand{\FF}{{\mathcal F}}

\newcommand{\QQ}{{\mathcal Q}}

\newcommand{\WW}{{\mathcal W}}

\newcommand{\gotD}{{\mathfrak D}}

\newcommand{\Fullbox}{{\rule{2.0mm}{2.0mm}}}

\newcommand{\EP}{\hfill\Fullbox\vspace{0.2cm}}
\newcommand{\prova}{\noindent{\it Proof. }}
\newcommand{\io}{\infty}
\newcommand{\eps}{\varepsilon}

\newcommand{\m}{\mu}

\newcommand{\g}{\gamma}
\newcommand{\om}{\omega}

\newcommand{\s}{\sigma}

\newcommand{\der}{{\rm d}}

\def\deriv#1#2#3{
\ifnum\catcode`#3=12    % If #3 is a number:
   \ifnum#3=1
      \frac {{\rm d}#1} {{\rm d}#2}
   \else
      \frac {{\rm d}^{#3}#1} {{\rm d}#2^{#3}}
   \fi
\else   % # 3 is a letter
      \frac {{\rm d}^{#3}#1} {{\rm d}#2^{#3}}
\fi}

\def\qed{\hfill\raise1pt\hbox{\vrule height5pt width5pt depth0pt}}

\def\ins#1#2#3{\vbox to0pt{\kern-#2 \hbox{\kern#1 #3}\vss}\nointerlineskip}

\begin{document}

\title{Frequency locking in the
injection-locked frequency divider equation}
\author{
\bf Michele V. Bartuccelli$^\ast$,
Jonathan H.B. Deane$^\ast$,
Guido Gentile$^\dagger$ 
\vspace{2mm}
\\ \small 
$^\ast$Department of Mathematics,
University of Surrey, Guildford, GU2 7XH, UK.
\\ \small 
E-mails:
m.bartuccelli@surrey.ac.uk, j.deane@surrey.ac.uk \\  \small  
$^\dagger$Dipartimento di Matematica, Universit\`a di Roma Tre, Roma,
I-00146, Italy.
\\ \small
E-mail: gentile@mat.uniroma3.it 
}

\date{}

\maketitle

\begin{abstract}  
We consider a model for the injection-locked frequency divider,
and study analytically the locking onto rational multiples of the
driving frequency. We provide explicit formulae for the width
of the plateaux appearing in the devil's staircase structure of the
lockings, and in particular show that the largest plateaux
correspond to even integer values for the ratio of the frequency
of the driving signal to the frequency of the output signal.
Our results prove the experimental
and numerical results available in the literature.
\end{abstract}

%\newpage

%\tableofcontents

%\newpage

%%%%%%%%%%%%%%%%%%%%%%%%%%%%%%%%%%%%%%%%%%%%%%%%%%%%%%%%%%%%%%%%%%%%%%%%%
%%%%%%%%%%%%%%%%%%%%%%%%%%%%%%%%%%%%%%%%%%%%%%%%%%%%%%%%%%%%%%%%%%%%%%%%%
\zerarcounters
\section{Introduction}
\label{sec:1}
%%%%%%%%%%%%%%%%%%%%%%%%%%%%%%%%%%%%%%%%%%%%%%%%%%%%%%%%%%%%%%%%%%%%%%%%%
%%%%%%%%%%%%%%%%%%%%%%%%%%%%%%%%%%%%%%%%%%%%%%%%%%%%%%%%%%%%%%%%%%%%%%%%%

In \cite{OBK}, an electronic circuit known as the injection-locked
frequency divider is studied experimentally, and the devil's staircase
structure of the lockings is measured: when the ratio of the
frequency $\om$ of the driving signal to the frequency $\Omega$
of the output signal is plotted versus $\om$, plateaux are found
for rational values of the ratio.
In \cite{OBYK}, a model for the circuit is presented and numerically
investigated, and the results are shown to agree with the experiments.

In this paper, on the basis of the model introduced in \cite{OBYK},
we address the problem of explaining analytically the appearance
of the plateaux of the devil's staircase. We aim to understand why
the largest plateaux correspond to even integer values for the
frequency ratio and, more generally, how the widths of the plateaux
depend on the particular values of the ratio.

From a qualitative point of view, the mechanism of locking
can be illustrated as follows. For fixed driving frequency $\om$
one considers the Poincar\'e section at times $t=2\pi n/\omega$,
for integer $n$, and studies the dynamics on the attractor.
This leads to a map which behaves as a diffeomorphism on the
circle. Thus, based on the theory of such systems \cite{A},
one expects that for $\om/\Omega_{0}$ close
to a rational number one has locking. How close $\om$ has to be to
a rational multiple of $\Omega_{0}$ depends on $\mu$ and
on the multiple itself: in the $(\om,\mu)$ parameter plane
one has locking in wedge-shaped regions known as Arnold tongues.

However, all the discussion above is purely qualitative.
In particular, there remains the major problem of determining
the map to which one should apply the theory.
A quantitative constructive analysis is another matter, and requires
taking into account the fine details of the equation and the
explicit expression of the solution of the unperturbed equation:
we carry this out in this paper. Our analysis is based
on perturbation theory, which is implemented to all orders
and proved to be convergent. This approach is particularly suited
for quantitative estimates within any given accuracy
(for which it has to be possible to go to arbitrarily
high perturbation orders, and to control the truncation errors).
Furthermore, we think that a rigorous analysis \textit{ab initio},
without introducing uncontrolled simplifications or approximations,
can be of interest by itself. Indeed, although such simplifications
can capture the essential features of the problem and allow a
qualitative understanding of the physical phenomenon, it
nonetheless remains unclear in general how far a simplified model
can be expected to describe the original system faithfully.

The conclusions of our analysis can be summarised as follows.
The equation modelling the system can be viewed as
a perturbation of order $\mu$ of a particular differential equation.
In the absence of the perturbation, after a suitable
change of variables, the system can be cast in the form of
a Li\'enard equation $x''+h(x)\,x'+k(x)=0$.
Under suitable assumptions on $h$ and $k$,
this admits a globally attracting limit cycle.
Let $\Omega_{0}$ be the proper frequency of such a cycle, and
let us denote by $x_{0}(t)=X_{0}(\Omega_{0}t)$ the solution
of the equation corresponding to the limit cycle,
with the function $X_{0}$ being $2\pi$-periodic in its argument.
By also including the time direction, one can study the
dynamics in the three-dimensional extended phase space $(x,x',t)$,
in which the limit cycle generates a topological cylinder.
When the perturbation is switched on, the cylinder survives
as an invariant manifold, slightly deformed with respect to
the unperturbed case. This follows from general arguments
related to the centre manifold theorem \cite{CH}.
However the dynamics on the manifold strongly depends on the relation
between the proper frequency $\Omega_{0}$ and the frequency $\om$ of
the driving signal. If $\om/\Omega_{0}$ is irrational
and satisfies some Diophantine condition (such as $|\om\nu_{1}+
\Omega_{0}\nu_{2}| > \g(|\nu_{1}|+|\nu_{2}|+1)^{-\tau}$ for all
$(\nu_{1},\nu_{2})\in\ZZZ^{2}$ and some positive constants $\g,\tau$),
then one expects the output signal $x(t)$ to be
a quasi-periodic function with frequency vector
$\boldsymbol\om=(\om,\Omega_{0})$, so that one
has $x(t)=X(\om t,\Omega_{0}t)=X_{0}(\Omega_{0}t)+O(\mu)$,
where $X$ is a $2\pi$-periodic function of both its arguments.
In this case we say that the output frequency $\Omega$
equals $\Omega_{0}$ (of course, this is slightly improper
terminology because $\Omega_{0}$ is only the frequency of the
leading contribution to the output signal, and the latter
is not even periodic).
On the other hand, if $\om/\Omega_{0}$ is close to a rational
number $p/q$ (resonance), then $x(t)$ is periodic with frequency
$\Omega=p\omega/q$ (locking): hence the frequency $\Omega$ of the
output signal differs from $\Omega_{0}$ --- even if it remains close
to it ---, because it is locked to the driving frequency $\om$.
Thus, if one plots the ratio $\omega/\Omega$ versus $\om$ one obtains
the devil's staircase structure depicted in Figures 4 to 9
of \cite{OBK}. The locked solutions can be obtained analytically
from the unperturbed periodic solutions by a mechanism similar
to the subharmonic bifurcations that we have studied in previous
papers \cite{BBDGG,BDG}. We stress, however, that, unlike
the cases studied in the latter references, here,
the unperturbed equation cannot be solved in closed form.
This will yield extra technical difficulties, because we shall have
to rely for our analysis on abstract symmetry properties of the
solution, without the possibility of using explicit expressions.

%%%%%%%%%%%%%%%%%%%%%%%%%%%%%%%%%%%%%%%%%%%%%%%%%%%%%%%%%%%%%%%%%%%%%%%%%
%%%%%%%%%%%%%%%%%%%%%%%%%%%%%%%%%%%%%%%%%%%%%%%%%%%%%%%%%%%%%%%%%%%%%%%%%
\zerarcounters
\section{Model for the injection-locked frequency divider}
\label{sec:2}
%%%%%%%%%%%%%%%%%%%%%%%%%%%%%%%%%%%%%%%%%%%%%%%%%%%%%%%%%%%%%%%%%%%%%%%%%
%%%%%%%%%%%%%%%%%%%%%%%%%%%%%%%%%%%%%%%%%%%%%%%%%%%%%%%%%%%%%%%%%%%%%%%%%

We consider the system of ordinary differential equations
\begin{equation}
C \frac{\der V_{C}}{\der t} = I_{L} + f(V_{C},t) , \qquad
L \frac{\der I_{L}}{\der t} = - R I_{L} - V_{C} , 
\label{eq:2.1}
\end{equation}
where $L,C,R>0$ are parameters, $V_{C}$ and $I_{L}$, the state variables, are
the capacitor voltage and the inductor current, respectively, and
\begin{equation}
f(V_{C},t)= \left( A + B \sin \Omega t \right)
V_{C} \left( 1- \left( V_{C}/V_{DD}\right)^{2} \right) ,
\qquad V_{DD},A>0 , \qquad B\in\RRR, \nonumber
\end{equation}
is the (cubic approximation of the)
driving point characteristic of the nonlinear resistor.
The model (\ref{eq:2.1}) was introduced in \cite{OBYK} as a simplified
description of the injection-locked frequency divider.

By introducing the new variables $u:=V_{C}/V_{DD}$ and
$v:=RI_{L}/V_{DD}$ and rescaling time $t\to R t/L$, 
(\ref{eq:2.1}) becomes
\begin{equation}
u' = \alpha v + \Phi(t) \, u \left( 1 - u^{2} \right) ,
\qquad v' = - u - v ,
\label{eq:2.2}
\end{equation}
where the prime denotes derivative with respect to time $t$,
and we have set $\alpha=L/R^{2}C$, $\beta=LA/RC$,
$\mu=LB/RC$, and $\Phi(t)=\beta + \mu \sin \omega t$,
with $\omega =\Omega L/R$.

From now on, we shall consider the system (\ref{eq:2.2}),
with $\alpha>\beta>1$, and $\mu,\omega\in\RRR$.
By setting $\sigma=u+v$ we obtain
\begin{equation}
\sigma' = \left(\alpha-1\right) \sigma +
\left( \Phi(t)-\alpha \right) u - \Phi(t) \,  u^{3} , \qquad
u' = \alpha \, \sigma +
\left(\Phi(t)-\alpha \right) u - \Phi(t) \, u^{3} ,
\label{eq:2.3}
\end{equation}
which gives $u'' + \left[ 1 - \Phi(t) + 3 \Phi(t) u^{2} \right] u' +
\left[ \left( \alpha - \Phi(t) \right) u + \Phi(t)\,u^{3} \right]  
+ \Phi'(t) \left( u^{3} - u \right) = 0$, that is
\begin{equation}
u'' + \left( 1 - \beta + 3 \beta u^{2} \right) u' +
\left[ \left( \alpha - \beta \right) u + \beta\,u^{3} \right]  
+ \mu \Psi(u,u',t) = 0 ,
\label{eq:2.4}
\end{equation}
with
\begin{equation}
\Psi(u,u',t)= \left[ u' \left( 3 u^{2} - 1 \right) \sin \omega t +
\left( u^{3} - u \right) \sin \omega t + \omega \left( u^{3} - u \right)
\cos \omega t \right] .
\label{eq:2.5}
\end{equation}
For $\m=0$, (\ref{eq:2.4}) reduces to
$u'' + ( 1 - \beta + 3 \beta u^{2}) \, u' + [(
\alpha - \beta ) \, u + \beta\,u^{3}] = 0$,
which can be written as a Li\'enard equation
\begin{equation}
u'' + u' \, h(u) + k(u) = 0 , 
\label{eq:2.6}
\end{equation}
with
\begin{equation}
h(u) = 1 - \beta +  3 \beta u^{2} ,
\qquad k(u) = \left( \alpha - \beta \right) u + \beta u^{3} =
u \left( \alpha - \beta + \beta u^{2} \right) .
\label{eq:2.7}
\end{equation}
For (\ref{eq:2.6}) to have a unique limit cycle \cite{C,H}, we require that
\begin{equation}
1-\beta < 0 , \qquad \alpha-\beta>0  \qquad
\Longrightarrow \qquad \alpha > \beta > 1 , \nonumber
\end{equation}
which motivates our assumptions on the parameters $\alpha$ and $\beta$.

\begin{figure}[htbp]
\centering
\vspace{0.3cm}
\includegraphics*[angle=0,width=2.5in]{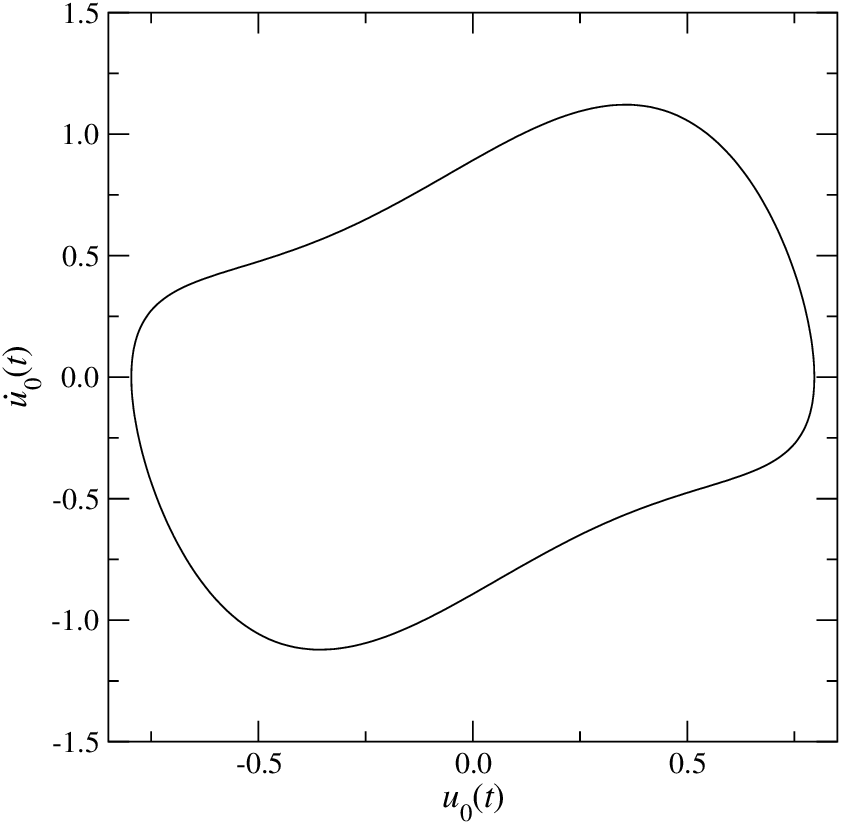}
\vspace{-0.3cm}
\caption{The limit cycle for $\alpha = 2.5$, $\beta = 2.0$
and $\mu=0$. The proper frequency is $\Omega_{0} \approx 1.1434$.}
\label{fig:1}
\end{figure}

Consider the system described by the equation (\ref{eq:2.6}),
with the functions $h(u)$ and $k(u)$ given by (\ref{eq:2.7})
with $\alpha > \beta > 1$.
Such a system admits one and only one limit cycle encircling
the origin \cite{H}; cf. Figure \ref{fig:1}. Let $T_{0}$ be the
period of the solution $u_{0}(t)$ running on such a cycle.
Denote by $\Omega_{0}=2\pi/T_{0}$ the corresponding frequency:
$\Omega_{0}$ will be called the \textit{proper frequency}
of the system. Note that $\Omega_{0}$ depends only on the
parameters $\alpha$ and $\beta$.

%%%%%%%%%%%%%%%%%%%%%%%%%%%%%%%%%%%%%%%%%%%%%%%%%%%%%%%%%%%%%%%%%%%%%%%%%
% FIGURE 1
%%%%%%%%%%%%%%%%%%%%%%%%%%%%%%%%%%%%%%%%%%%%%%%%%%%%%%%%%%%%%%%%%%%%%%%%%
%%%%%%%%%%%%%%%%%%%%%%%%%%%%%%%%%%%%%%%%%%%%%%%%%%%%%%%%%%%%%%%%%%%%%%%%%

The solution $u_{0}(t)$ is unique up to time translation.
Fix the time origin so that $u'_{0}(0)=0$, $u_{0}(0)>0$.
Note that fixing the origin of time in such a way
that $u'_{0}(0)=0$ compels us to shift by some $t_{0}$
the time in the argument of the driving term in (\ref{eq:2.5}),
i.e. $\Psi(u,u',t)$ must be replaced with $\Psi(u,u',t+t_{0})$;
cf. the analogous discussion in \cite{GBD3}.

%%%%%%%%%%%%%%%%%%%%%%%%%%%%%%%%%%%%%%%%%%%%%%%%%%%%%%%%%%%%%%%%%%%%%%%%%
% LEMMA 1
%%%%%%%%%%%%%%%%%%%%%%%%%%%%%%%%%%%%%%%%%%%%%%%%%%%%%%%%%%%%%%%%%%%%%%%%%
\begin{lemma} \label{lem:3.1}
The Fourier expansion of $u_{0}(t)$ contains only the odd harmonics, i.e.
\begin{equation}
u_{0}(t) = \sum_{\nu\in\ZZZ} {\rm e}^{i \Omega_{0} \nu t}
u_{0,\nu} = \sum_{\substack{\nu\in\ZZZ\\\nu \; {\rm odd}}}
{\rm e}^{i \Omega_{0} \nu t} u_{0,\nu} .
\nonumber
\end{equation}
\end{lemma}
%%%%%%%%%%%%%%%%%%%%%%%%%%%%%%%%%%%%%%%%%%%%%%%%%%%%%%%%%%%%%%%%%%%%%%%%%

%%%%%%%%%%%%%%%%%%%%%%%%%%%%%%%%%%%%%%%%%%%%%%%%%%%%%%%%%%%%%%%%%%%%%%%%%
\prova The symmetry properties of (\ref{eq:2.6}),
more precisely the fact that $h(-u) = h(u)$ and $k(-u) = -k(u)$,
ensure that the periodic solution $u_{0}(t)$ satisfies the property
\begin{equation}
u_{0}(t + T_{0}/2) = -u_{0}(t) ,
\label{eq:3.3} 
\end{equation}
and in turn this implies the result (compare
the proof of Lemma 3.2 in \cite{BDG}).\EP
%%%%%%%%%%%%%%%%%%%%%%%%%%%%%%%%%%%%%%%%%%%%%%%%%%%%%%%%%%%%%%%%%%%%%%%%%

%%%%%%%%%%%%%%%%%%%%%%%%%%%%%%%%%%%%%%%%%%%%%%%%%%%%%%%%%%%%%%%%%%%%%%%%%
% LEMMA 2
%%%%%%%%%%%%%%%%%%%%%%%%%%%%%%%%%%%%%%%%%%%%%%%%%%%%%%%%%%%%%%%%%%%%%%%%%
\begin{lemma} \label{lem:3.2}
One has $\displaystyle \int_{0}^{T_{0}} {\rm d}t \, h(u_{0}(t)) > 0$.
\end{lemma}
%%%%%%%%%%%%%%%%%%%%%%%%%%%%%%%%%%%%%%%%%%%%%%%%%%%%%%%%%%%%%%%%%%%%%%%%%

%%%%%%%%%%%%%%%%%%%%%%%%%%%%%%%%%%%%%%%%%%%%%%%%%%%%%%%%%%%%%%%%%%%%%%%%%
\prova For a proof see \cite{C}.\EP
%%%%%%%%%%%%%%%%%%%%%%%%%%%%%%%%%%%%%%%%%%%%%%%%%%%%%%%%%%%%%%%%%%%%%%%%%
\vskip.2truecm
%%%%%%%%%%%%%%%%%%%%%%%%%%%%%%%%%%%%%%%%%%%%%%%%%%%%%%%%%%%%%%%%%%%%%%%%%

\begin{figure}[htbp]
\vspace{0.3cm}
\centering
\subfigure[]{
    	\includegraphics*[angle=0,width=2.5in]{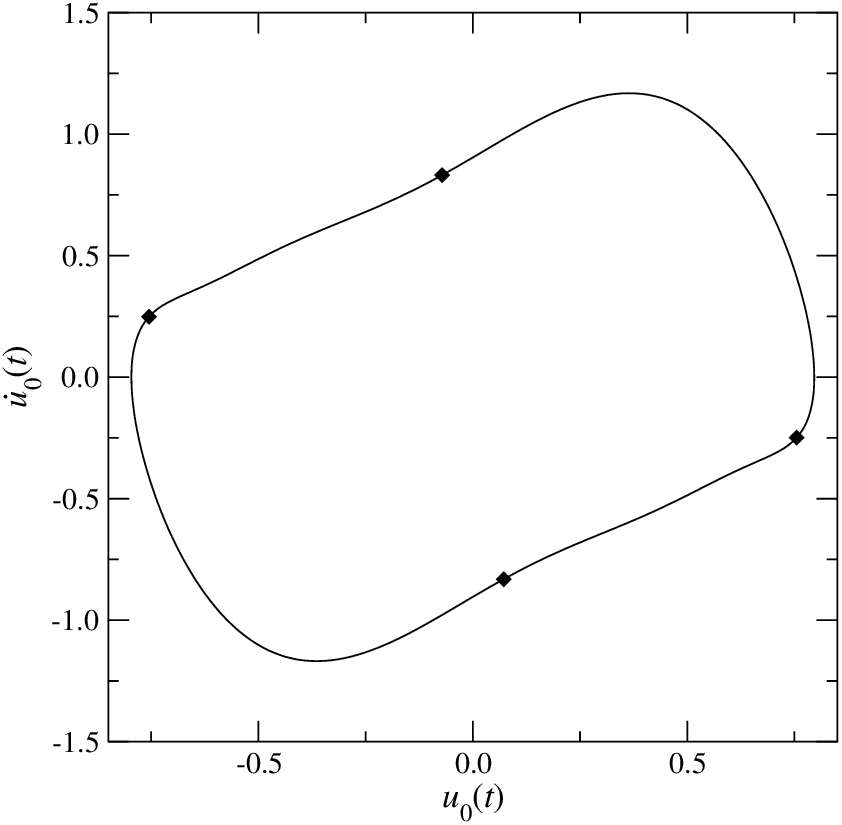}}
\hspace{0.5in}
\subfigure[]{
    	\includegraphics*[angle=0,width=2.5in]{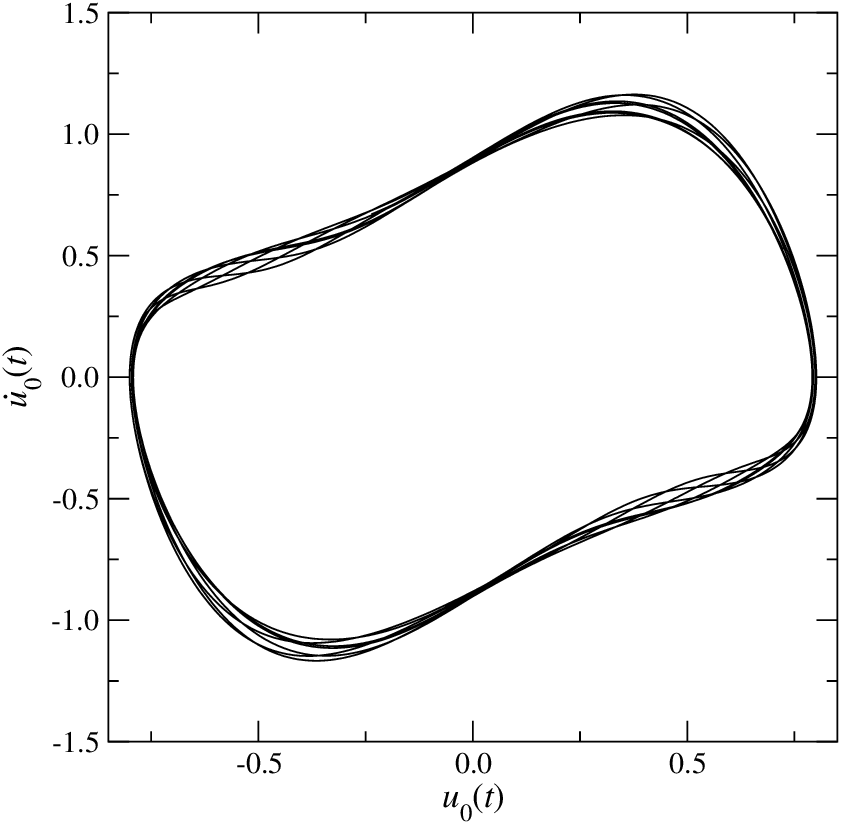}}
\vspace{-0.3cm}
\caption{Examples of attractors for $\alpha = 2.5$, $\beta = 2.0$
and $\mu=0.1$. For $|\omega-4\Omega_{0}| \le 0.03$
the motion is periodic: in (a), $\omega=4\Omega_{0} + 0.02$.
The black diamonds mark the four points where $\sin\omega t$ is zero
and positive-going.
For $\omega=4\Omega_{0} + 0.2$ the motion is quasi-periodic (b).
Recall that $\Omega_{0} \approx 1.1434$.}
\label{fig:2}
\end{figure}

Moreover the limit cycle is a global attractor \cite{H,SC},
and it is uniformly hyperbolic \cite{Z,C}. Hence the cylinder
it generates in the extended phase space persists, slightly deformed,
as a global attractor for small perturbations \cite{L2,CH,HPS,V}.
This also means that the system described by the equation (\ref{eq:2.4}),
at least for small values of $\mu$, has one and only one attractor,
and the latter attracts the whole phase space.
However, the persistence of the attractor does not tell us
whether the dynamics on the attractor is periodic or
quasi-periodic; cf. \cite{CLB} for an analogous discussion.
In particular it does not imply that for $\Omega_{0}/\omega$
close to a resonance the dynamics remains periodic;
cf. Figure \ref{fig:2}.

%%%%%%%%%%%%%%%%%%%%%%%%%%%%%%%%%%%%%%%%%%%%%%%%%%%%%%%%%%%%%%%%%%%%%%%%%
% FIGURE 2
%%%%%%%%%%%%%%%%%%%%%%%%%%%%%%%%%%%%%%%%%%%%%%%%%%%%%%%%%%%%%%%%%%%%%%%%%
%%%%%%%%%%%%%%%%%%%%%%%%%%%%%%%%%%%%%%%%%%%%%%%%%%%%%%%%%%%%%%%%%%%%%%%%%

We note that for $\Omega_{0}/\omega$ Diophantine the attractor is
expected to become quasi-periodic, with the dynamics
analytically conjugated to a Diophantine rotation with
rotation vector $(\Omega_{0},\omega)$. In principle, this can be
proved by KAM techniques \cite{CLB}, or with methods closer
to those used in this paper \cite{G1,G2,GG,GBD1,GBG}.

%%%%%%%%%%%%%%%%%%%%%%%%%%%%%%%%%%%%%%%%%%%%%%%%%%%%%%%%%%%%%%%%%%%%%%%%%
%%%%%%%%%%%%%%%%%%%%%%%%%%%%%%%%%%%%%%%%%%%%%%%%%%%%%%%%%%%%%%%%%%%%%%%%%
\zerarcounters
\section{Framework for studying frequency locking}
\label{sec:5}
%%%%%%%%%%%%%%%%%%%%%%%%%%%%%%%%%%%%%%%%%%%%%%%%%%%%%%%%%%%%%%%%%%%%%%%%%
%%%%%%%%%%%%%%%%%%%%%%%%%%%%%%%%%%%%%%%%%%%%%%%%%%%%%%%%%%%%%%%%%%%%%%%%%

Rescale time so that the driving term has period $2\pi$, hence
frequency $1$, by setting $\tau=\omega t$. Then, by denoting with
the dot the derivative with respect to rescaled time $\tau$,
(\ref{eq:2.4}) gives
\begin{equation} \ddot u + \frac{1}{\omega} 
\left( 1 - \beta + 3 \beta u^{2} \right) \dot u +
\frac{1}{\omega^{2}}
\left[ \left( \alpha - \beta \right) u + \beta\,u^{3} \right]  
+ \mu \bar \Psi(u,\dot u,\tau+\tau_{0}) = 0 ,
\label{eq:4.1}
\end{equation}
where $\tau_{0}=\omega t_{0}$, and we have defined
$\bar \Psi(u,\dot u,\tau) = [ \omega^{-1} \,
\dot u \, ( 3 u^{2} - 1 ) \, \sin \tau $ $+$
$\omega^{-2}( u^{3} - u ) \, \sin \tau +
\omega^{-1} ( u^{3} - u ) \, \cos \tau ]$.
For $\mu=0$ one has $\ddot u + \omega^{-1} 
( 1 - \beta + 3 \beta u^{2} ) \, \dot u +
\omega^{-2} [ ( \alpha - \beta ) \,
u + \beta\,u^{3} ]  = 0$, which can be written as
\begin{equation}
\ddot u + \frac{1}{\omega} \, h(u) \, \dot u +
\frac{1}{\omega^{2}} \, k(u) = 0 ,
\label{eq:4.3}
\end{equation}
which is of the form (\ref{eq:2.6}) up to the rescaling of time.
As an effect of the time rescaling, the frequency of the limit cycle
for the system (\ref{eq:4.3}) depends on $\omega$, as it is given
by $\overline\Omega_{0}=\Omega_{0}/\omega$.

%%%%%%%%%%%%%%%%%%%%%%%%%%%%%%%%%%%%%%%%%%%%%%%%%%%%%%%%%%%%%%%%%%%%%%%%%
% REMARK 1
%%%%%%%%%%%%%%%%%%%%%%%%%%%%%%%%%%%%%%%%%%%%%%%%%%%%%%%%%%%%%%%%%%%%%%%%%
\begin{rmk} \label{rmk:1}
As the solution $u_{0}(\tau)$ is analytic in $\tau$, the property
$\dot u_{0}(0)=0$ means that we can write
$\dot u_{0}(\tau) = r_{1} \tau + O(\tau^{2})$, and hence
$u_{0}(\tau) = r_{0} + r_{1}\tau^{2}/2 +O(\tau^{3})$,
with $r_{0}=u_{0}(0)$ and $r_{1}=\ddot u_{0}(0)$.
\end{rmk}
%%%%%%%%%%%%%%%%%%%%%%%%%%%%%%%%%%%%%%%%%%%%%%%%%%%%%%%%%%%%%%%%%%%%%%%%%

We want to show that if the frequency $\omega$ of the driving term
is close to a rational multiple of the unperturbed proper frequency
$\Omega_{0}$ of the system, that is $\omega \approx p \Omega_{0}/q$
for some $p,q\in\NNN$ relatively prime, then the frequency $\Omega$
of the solution exactly equals $q\omega/p$, that is $\omega/\Omega=p/q$.
Such a phenomenon is known as \textit{frequency locking}:
the system is said to be locked into the \textit{resonance} $p\!:\!q$.

Let $\rho=p/q\in\QQQ$. For $\mu=0$, for any frequency $\omega$ of the
driving term the proper frequency is $\Omega=\Omega_{0}$ --- the
system is decoupled from the perturbation ---, so that if we fix
$\omega=\rho\Omega_{0}$ we obtain $\Omega=\Omega_{0}=\omega/\rho$.
In terms of the rescaled variables, for which $\omega$ is replaced with
$\overline\omega=1$, the proper frequency
becomes $\overline\Omega_{0}=1/\rho$.
For $\omega$ close to $\rho\Omega_{0}$ write
\begin{equation}
\frac{1}{\omega} = \frac{1}{\rho\Omega_{0}} + \eps(\mu) ,
\label{eq:5.1}
\end{equation}
with $\eps(\mu)$ such that $\eps(\mu)\to0$ as $\mu\to0$.

We look for periodic solutions for the full system (\ref{eq:2.4}),
hence for solutions with period $T=2\pi p/\omega$ (i.e the
least common multiple of both $2\pi/\omega$ and $2\pi p/\omega q$).
In terms of the rescaled time $\tau$, the solution will have period
$2\pi p$, hence frequency $1/p$. For $\mu=0$ the system (\ref{eq:4.1})
reduces to
\begin{equation}
H_{0}(u,\dot u,\ddot u) := \ddot u + f(u) \, \dot u + g(u) = 0 ,
\label{eq:5.2}
\end{equation}
with
\begin{eqnarray}
f(u) = \frac{1}{\rho\Omega_{0}} \, h(u) , \qquad
g(u) = \frac{1}{\rho^{2}\Omega_{0}^{2}} \, k(u) , \nonumber
\end{eqnarray}
which admits the periodic solution $u_{0}(\tau)$ such that
$u_{0}(0)>0$, $\dot u_{0}(0)=0$ and $u_{0}(\tau+2\pi\rho)=u_{0}(\tau)$.
In other words the frequency of the limit cycle is
$1/\rho=q/p$ and the period is $2\pi p/q$, i.e. $u_{0}(\tau) =
U(\tau/\rho)$, with the function $U$ being $2\pi$-periodic.

For $\mu\neq0$ we write
\begin{equation}
\eps(\mu) = \eps_{1} \mu + \eps_{2} \mu^{2} + \ldots =
\sum_{k=1}^{\io} \eps_{k} \mu^{k} ,
\label{eq:5.5}
\end{equation}
and, by inserting (\ref{eq:5.1}) and (\ref{eq:5.5})
into (\ref{eq:4.1}), we obtain the equation
\begin{equation}
H(u,\dot u,\ddot u,\mu) := H_{0}(u,\dot u,\ddot u) +
\sum_{k=1}^{\io} \mu^{k} H_{k} (u,\dot u,\tau+\tau_{0}) = 0 , 
\label{eq:5.6}
\end{equation}
where
\begin{eqnarray}
\null \hskip-.7truecm H_{1}(u,\dot u,\tau) & = &
\eps_{1} \left( 1 - \beta + 3 \beta u^{2} \right) \dot u +
\frac{2\eps_{1} }{\rho\Omega_{0}} \left[ \left( \alpha - \beta \right)
u + \beta\,u^{3} \right] \nonumber \\
\null \hskip-.7truecm & + & 
\frac{1}{\rho\Omega_{0}} \, \dot u \left(
3 u^{2} - 1 \right) \sin \tau +
\frac{1}{\rho^{2}\Omega_{0}^{2}} \left( u^{3} - u \right) \sin \tau
+ \frac{1}{\rho\Omega_{0}} \left( u^{3} - u \right) \cos \tau ,
\label{eq:5.7} \\
\null \hskip-.7truecm H_{2}(u,\dot u,\tau) & = &
\eps_{2} \left( 1 - \beta + 3 \beta u^{2} \right) \dot u +
\left( \frac{2\eps_{2}}{\rho\Omega_{0}} + \eps_{1}^{2} \right)
\left[ \left( \alpha - \beta \right) u +
\beta\,u^{3} \right] \nonumber \\
\null \hskip-.7truecm & + &
\eps_{1} \dot u \left( 3 u^{2} - 1 \right) \sin \tau +
\frac{2\eps_{1}}{\rho\Omega_{0}} \left( u^{3} - u \right) \sin \tau +
\eps_{1} \left( u^{3} - u \right) \cos \tau ,
\label{eq:5.8}
\end{eqnarray}
and so on. The shifting of time by $\tau_{0}=\omega t_{0}$ in the
driving term is due to the choice of the origin of time made
according to Section \ref{sec:2}.

In the following sections we shall prove that for $\mu$ small enough
it is possible to choose $\eps(\mu)$ as a function of $t_{0}$,
in such a way that there exists a periodic solution of (\ref{eq:5.6})
with period $2\pi p$, i.e. with frequency $1/p$. When projected onto
the $(u,\dot u)$ plane, such a solution is close enough to the
unperturbed limit cycle (cf. for instance Figure \ref{fig:2}):
the difference between them is of order $\mu$.

%%%%%%%%%%%%%%%%%%%%%%%%%%%%%%%%%%%%%%%%%%%%%%%%%%%%%%%%%%%%%%%%%%%%%%%%%
%%%%%%%%%%%%%%%%%%%%%%%%%%%%%%%%%%%%%%%%%%%%%%%%%%%%%%%%%%%%%%%%%%%%%%%%%
\zerarcounters
\section{The linearised equation}
\label{sec:6}
%%%%%%%%%%%%%%%%%%%%%%%%%%%%%%%%%%%%%%%%%%%%%%%%%%%%%%%%%%%%%%%%%%%%%%%%%
%%%%%%%%%%%%%%%%%%%%%%%%%%%%%%%%%%%%%%%%%%%%%%%%%%%%%%%%%%%%%%%%%%%%%%%%%

Write the unperturbed system (\ref{eq:5.2}) as
\begin{equation}
\dot u = v , \qquad
\dot v = G(u,v) ,
\label{eq:6.1}
\end{equation}
with $G(u,v) = - (\rho\Omega_{0})^{-1} 
( 1 - \beta + 3 \beta u^{2}) \, v - (\rho^{2}\Omega_{0})^{-2}
[ ( \alpha - \beta ) \, u + \beta\,u^{3}]$.
Let $(u_{0}(\tau),v_{0}(\tau))$ be the $2\pi \rho$-periodic solution
of (\ref{eq:6.1}), which is uniquely determined by the conditions
\begin{equation} \dot u_{0}(0) = 0 , \qquad u_{0}(0) > 0 .
\label{eq:6.3}
\end{equation}
The periodicity properties of $u_{0}(\tau)$ allow us to write
\begin{equation}
u_{0}(\tau) = \sum_{\substack{{\nu\in\ZZZ} \\ \nu \; {\rm odd}}}
{\rm e}^{i\nu \tau /\rho} u_{0,\nu} , \qquad \rho = \frac{p}{q} ,
\label{eq:6.4}
\end{equation}
as follows from Lemma \ref{lem:3.1}. Denote by
\begin{equation}
W(\tau) = \left( \begin{matrix}
w_{11}(\tau) & w_{12}(\tau) \\
w_{21}(\tau) & w_{22}(\tau) \end{matrix} \right) 
\label{eq:6.5}
\end{equation}
the Wronskian matrix of the system (\ref{eq:6.1}), that is
the solution of the matrix equation
\begin{equation}
\begin{cases}
\dot W(\tau) = M(\tau) \, W(\tau) , & \\ W(0) = \one , \end{cases} \qquad
M(\tau) = \left( \begin{matrix}
0 & 1 \\ G_{u}(u_{0}(\tau),v_{0}(\tau)) & G_{v}(u_{0}(\tau),v_{0}(\tau))
\end{matrix} \right) ,
\label{eq:6.6}
\end{equation}
where $G_{u}$ and $G_{v}$ denote derivatives with respect to
$u$ and $v$ of $G$, and $w_{21}(\tau)=\dot w_{11}(\tau)$,
$w_{22}(\tau)=\dot w_{12}(\tau)$.

%%%%%%%%%%%%%%%%%%%%%%%%%%%%%%%%%%%%%%%%%%%%%%%%%%%%%%%%%%%%%%%%%%%%%%%%%
% LEMMA 3
%%%%%%%%%%%%%%%%%%%%%%%%%%%%%%%%%%%%%%%%%%%%%%%%%%%%%%%%%%%%%%%%%%%%%%%%%
\begin{lemma} \label{lem:6.1}
In (\ref{eq:6.5}) one can set
\begin{equation}
w_{12}(\tau) := c_{2} \dot u_{0}(\tau) , \qquad
w_{11}(\tau) := c_{1} \dot u_{0}(\tau) \int_{\bar \tau}^{\tau}
{\rm d} \tau' \frac{{\rm e}^{-F(\tau')}}{\dot u_{0}^{2}(\tau')} , 
\label{eq:6.7}
\end{equation}
where $F(\tau)$ is defined as
\begin{equation}
F(\tau) := \int_{0}^{\tau} {\rm d}\tau' \, f(u_{0}(\tau')) , \qquad
f(u) := \frac{1}{\rho\Omega_{0}} 
\left( 1 - \beta + 3 \beta u^{2} \right) =
\frac{1}{\rho\Omega_{0}} \, h(u) ,
\label{eq:6.8}
\end{equation}
the constant $\bar \tau \in (0,\pi\rho)$ is chosen so that
$\dot w_{11}(0)=0$, and the constants $c_{1}$ and $c_{2}$ 
are such that $w_{11}(0)=w_{22}(0)=1$.
\end{lemma}
%%%%%%%%%%%%%%%%%%%%%%%%%%%%%%%%%%%%%%%%%%%%%%%%%%%%%%%%%%%%%%%%%%%%%%%%%

%%%%%%%%%%%%%%%%%%%%%%%%%%%%%%%%%%%%%%%%%%%%%%%%%%%%%%%%%%%%%%%%%%%%%%%%%
\prova It can immediately be checked that $(w_{12}(\tau),w_{22}(\tau))$,
with $w_{12}(\tau)$ defined as in (\ref{eq:6.7}) and $w_{22}(\tau)=
\dot w_{12}(\tau)$,
solves the linearised equation to (\ref{eq:6.1}). Then a second
independent solution is of the form $(w_{11}(\tau),w_{21}(\tau))$,
with $w_{11}(\tau)$ given by (\ref{eq:6.7}) and
$w_{21}(\tau)=\dot w_{11}(\tau)$; cf. \cite{I}, p. 122.
In Appendix \ref{app:A} we show that it is possible to choose
$\bar \tau \in (0,\pi\rho)$ in such a way that $\dot w_{11}(0)=0$.
The constants $c_{1}$ and $c_{2}$ are chosen so that $W(0)=\one$.\EP
%%%%%%%%%%%%%%%%%%%%%%%%%%%%%%%%%%%%%%%%%%%%%%%%%%%%%%%%%%%%%%%%%%%%%%%%%

%%%%%%%%%%%%%%%%%%%%%%%%%%%%%%%%%%%%%%%%%%%%%%%%%%%%%%%%%%%%%%%%%%%%%%%%%
% REMARK 2
%%%%%%%%%%%%%%%%%%%%%%%%%%%%%%%%%%%%%%%%%%%%%%%%%%%%%%%%%%%%%%%%%%%%%%%%%
\begin{rmk} \label{rmk:2}
With the notations of Remark \ref{rmk:1} one has $c_{2}=1/r_{1}$
and $c_{1}=-r_{1}$, so that $c_{1}c_{2}+1=0$.
\end{rmk}
%%%%%%%%%%%%%%%%%%%%%%%%%%%%%%%%%%%%%%%%%%%%%%%%%%%%%%%%%%%%%%%%%%%%%%%%%

Note that $\dot u_{0}(0)=0$, so that in (\ref{eq:6.7})
the function $w_{11}(\tau)$ at $\tau=0$ is defined as the limit
\begin{equation}
\lim_{\tau\to 0} c_{1} \dot u_{0}(\tau)
\int_{\bar\tau}^{\tau} {\rm d} \tau'
\frac{{\rm e}^{-F(\tau')}}{\dot u_{0}^{2}(\tau')} , 
\label{eq:6.9}
\end{equation}
which is well defined; cf. Appendix \ref{app:A}.
The same argument applies for $\tau=\pi\rho$, where
$\dot u_{0}(\pi\rho)= \dot u_{0}(0)=0$ --- by (\ref{eq:3.3}),
with the half-period $T_{0}/2$ becoming $\pi\rho$ in terms of
the rescaled variable.

For any periodic function $G$ we denote its average by
$\langle G \rangle$ and set $\tilde G = G - \langle G \rangle$.
Then $f_{0}:=\langle f\circ u_{0} \rangle > 0$
(cf. Lemma \ref{lem:3.2}), so that we can write
\begin{equation}
{\rm e}^{-F(\tau)} = {\rm e}^{-f_{0}\tau - \tilde F(\tau)} ,
\qquad \tilde F(\tau) = \int_{0}^{\tau} {\rm d}\tau'
\left( f(u_{0}(\tau'))- f_{0} \right) , 
\label{eq:6.10}
\end{equation}
where $\tilde F(\tau)$, and hence ${\rm e}^{-\tilde F(\tau)}$,
are well defined $2\pi\rho\,$-periodic functions.

%%%%%%%%%%%%%%%%%%%%%%%%%%%%%%%%%%%%%%%%%%%%%%%%%%%%%%%%%%%%%%%%%%%%%%%%%
% LEMMA 4
%%%%%%%%%%%%%%%%%%%%%%%%%%%%%%%%%%%%%%%%%%%%%%%%%%%%%%%%%%%%%%%%%%%%%%%%%
\begin{lemma} \label{lem:6.2}
Given any periodic function $P(\tau)$ and any real constant $C\neq0$
there exists a periodic function $Q(\tau)$, with the same period
as $P(\tau)$, and a constant $D$ such that
\begin{equation}
\int_{0}^{\tau} {\rm d}\tau' \, {\rm e}^{C\tau'} P(\tau') =
D + {\rm e}^{C \tau} Q(\tau) . \nonumber
\end{equation}
One has $D=-Q(0)$.
\end{lemma}
%%%%%%%%%%%%%%%%%%%%%%%%%%%%%%%%%%%%%%%%%%%%%%%%%%%%%%%%%%%%%%%%%%%%%%%%%

%%%%%%%%%%%%%%%%%%%%%%%%%%%%%%%%%%%%%%%%%%%%%%%%%%%%%%%%%%%%%%%%%%%%%%%%%
\prova Let $P(\tau)$ be a periodic function of period $T$. Write
\begin{equation}
P(\tau) = \sum_{\nu\in\ZZZ} {\rm e}^{i \omega \nu \tau} P_{\nu} ,
\label{eq:6.11} \end{equation}
where $\omega=2\pi/T$. Then one has
\begin{equation}
\int_{0}^{\tau} {\rm d}\tau' \, {\rm e}^{C\tau'} P(\tau') =
\sum_{\nu\in\ZZZ} P_{\nu} \int_{0}^{\tau} {\rm d}\tau' \,
{\rm e}^{i \omega \nu \tau' + C\tau'} =
\sum_{\nu\in\ZZZ} P_{\nu} \frac{ {\rm e}^{i \omega \nu \tau +
C \tau} - 1}{C+i\omega\nu} , \nonumber
\end{equation}
so that, by setting
\begin{equation}
Q(\tau) := \sum_{\nu\in\ZZZ}
\frac{P_{\nu}}{C+i\omega\nu} {\rm e}^{i \omega \nu \tau} , \qquad
D := - \sum_{\nu\in\ZZZ} \frac{P_{\nu}}{C+i\omega\nu} ,
\label{eq:6.12} \end{equation}
the assertion follows.\EP
%%%%%%%%%%%%%%%%%%%%%%%%%%%%%%%%%%%%%%%%%%%%%%%%%%%%%%%%%%%%%%%%%%%%%%%%%

%%%%%%%%%%%%%%%%%%%%%%%%%%%%%%%%%%%%%%%%%%%%%%%%%%%%%%%%%%%%%%%%%%%%%%%%%
% LEMMA 5
%%%%%%%%%%%%%%%%%%%%%%%%%%%%%%%%%%%%%%%%%%%%%%%%%%%%%%%%%%%%%%%%%%%%%%%%%
\begin{lemma} \label{lem:6.3}
There exist two $2\pi\rho$-periodic functions
$a(\tau)$ and $b(\tau)$ such that
\begin{equation}
w_{11}(\tau) = a(\tau) + {\rm e}^{-f_{0}\tau} \, b(\tau) , \qquad
w_{12}(\tau) = c \, a(\tau) ,
\label{eq:6.13}
\end{equation}
for a suitable constant $c$.
\end{lemma}
%%%%%%%%%%%%%%%%%%%%%%%%%%%%%%%%%%%%%%%%%%%%%%%%%%%%%%%%%%%%%%%%%%%%%%%%%

%%%%%%%%%%%%%%%%%%%%%%%%%%%%%%%%%%%%%%%%%%%%%%%%%%%%%%%%%%%%%%%%%%%%%%%%%
\prova We cannot directly apply Lemma \ref{lem:6.2} because
the function ${\rm e}^{-F(\tau)}/\dot u_{0}^{2}(\tau)$ appearing
in (\ref{eq:6.7}) is singular. However we can proceed as follows.
We write
$$ \frac{1}{\dot u_{0}^{2}(\tau)} = \lim_{\eta\to0}
\frac{1}{\dot u_{0}^{2}(\tau)+\eta} , $$
so that the new integrand is smooth and it is given by
${\rm e}^{-f_{0}\tau}$ times a $2\pi\rho$-periodic function.
Hence, as long as $0<\tau<\pi\rho$, the integrand is
bounded uniformly in $\eta$, and we can apply
Lebesgue's dominated convergence theorem, to write
$$ w_{11}(\tau) = c_{1} \dot u_{0}(\tau) \lim_{\eta\to0}
\int_{\bar\tau}^{\tau} {\rm d}\tau' \,
\frac{{\rm e}^{ct'}}{\dot u_{0}^{2}(\tau')+\eta} . $$
Then Lemma \ref{lem:6.2} gives
$$ w_{11}(\tau) = c_{1} \dot u_{0}(\tau) 
\lim_{\eta\to0} \left( {\rm e}^{c\tau} P(\tau,\eta) -
{\rm e}^{c \bar\tau} P(\bar\tau,\eta) \right) =
c_{1} \lim_{\eta\to0} \left( \dot u_{0}(\tau)
{\rm e}^{c\tau} P(\tau,\eta) \right) -
\dot u_{0}(\tau) {\rm e}^{c \bar\tau} P(\bar\tau) , $$
where the function $P(\tau,\eta)$ is $2\pi\rho\,$-periodic in $\tau$
and $P(\bar\tau)=\lim_{\eta\to0}P(\bar\tau,\eta)$ is well defined.
Note that  $-\dot u_{0}(\tau){\rm e}^{c \bar\tau} P(\bar\tau)$ gives
the function $a(\tau)$ in (\ref{eq:6.13}).
On the other hand, the function $w_{11}(\tau)$ is also well defined,
so that we can conclude that
$\lim_{\eta\to0} \left( \dot u_{0}(\tau)
{\rm e}^{c\tau} P(\tau,\eta) \right) $
is well defined and smooth. As the function $\dot u_{0}(\tau)\,
P(\tau,\eta)$ is periodic for any $\eta$, the limit will also be
periodic, and this defines the function $b(\tau)$ of (\ref{eq:6.13}).

Comparing the expressions for $w_{11}(\tau)$ and $w_{12}(\tau)$ in
(\ref{eq:6.7}), proportionality between the function $w_{12}(\tau)$
and the periodic component of $w_{11}(\tau)$ also follows.\EP
%%%%%%%%%%%%%%%%%%%%%%%%%%%%%%%%%%%%%%%%%%%%%%%%%%%%%%%%%%%%%%%%%%%%%%%%%

%%%%%%%%%%%%%%%%%%%%%%%%%%%%%%%%%%%%%%%%%%%%%%%%%%%%%%%%%%%%%%%%%%%%%%%%%
% LEMMA 6
%%%%%%%%%%%%%%%%%%%%%%%%%%%%%%%%%%%%%%%%%%%%%%%%%%%%%%%%%%%%%%%%%%%%%%%%%
\begin{lemma} \label{lem:6.4}
The Fourier expansions of the functions $a(\tau)$ and $b(\tau)$
in (\ref{eq:6.13}) contain only the odd harmonics.
\end{lemma}
%%%%%%%%%%%%%%%%%%%%%%%%%%%%%%%%%%%%%%%%%%%%%%%%%%%%%%%%%%%%%%%%%%%%%%%%%

%%%%%%%%%%%%%%%%%%%%%%%%%%%%%%%%%%%%%%%%%%%%%%%%%%%%%%%%%%%%%%%%%%%%%%%%%
\prova Write $u_{0}(\tau)$ according to (\ref{eq:6.4}). Then
$w_{12}(\tau) = c \, a(\tau) = c \, \dot u_{0}(\tau)$, so that
the assertion follows trivially for $a(\tau)$.
Moreover the function ${\rm e}^{-F(\tau)}/\dot u_{0}^{2}(\tau)$
involves even powers of functions containing only odd harmonics,
so that it contains only even harmonics, and so does its integral
as appearing in the definition (\ref{eq:6.7}) of $w_{11}(\tau)$.
Hence, by Lemma \ref{lem:6.2} and Lemma \ref{lem:6.3},
also $b(\tau)$ in (\ref{eq:6.13}) contains only the odd harmonics.\EP
%%%%%%%%%%%%%%%%%%%%%%%%%%%%%%%%%%%%%%%%%%%%%%%%%%%%%%%%%%%%%%%%%%%%%%%%%
\vskip.2truecm
%%%%%%%%%%%%%%%%%%%%%%%%%%%%%%%%%%%%%%%%%%%%%%%%%%%%%%%%%%%%%%%%%%%%%%%%%

A straightforward calculation gives
$\det W(\tau) = - c_{1} c_{2} {\rm e}^{-F(\tau)} = {\rm e}^{-F(\tau)} $,
since $c_{1}c_{2}=-1$, so that
\begin{equation}
W^{-1}(\tau) = {\rm e}^{F(\tau)} \left( \begin{matrix}
w_{22}(\tau) & -w_{12}(\tau) \\
-w_{21}(\tau) & w_{11}(\tau) \end{matrix} \right) .
\label{eq:6.15}
\end{equation}
We want to develop perturbation theory for a $2\pi p$-periodic solution
which continues the solution running on the unperturbed limit cycle
when the perturbation is switched on. Therefore we write
\begin{equation}
u(\tau) = u_{0}(\tau) + \sum_{k=1}^{\io} \mu^{k} u_{k}(\tau) , \qquad
u_{k}(\tau) = \sum_{\nu\in\ZZZ} {\rm e}^{i\nu \tau/p} u_{k,\nu} ,
\label{eq:6.16}
\end{equation}
where $u_{0}(\tau)$ is the solution satisfying the conditions
(\ref{eq:6.3}). Inserting (\ref{eq:6.16}) into (\ref{eq:5.6})
and expanding everything in powers of $\mu$,
we obtain a sequence of recursive equations.
In Sections \ref{sec:7} and \ref{sec:8} we shall consider
in detail the first order. Higher order analysis
and the issue of convergence will be discussed in Section \ref{sec:10}.

%%%%%%%%%%%%%%%%%%%%%%%%%%%%%%%%%%%%%%%%%%%%%%%%%%%%%%%%%%%%%%%%%%%%%%%%%
%%%%%%%%%%%%%%%%%%%%%%%%%%%%%%%%%%%%%%%%%%%%%%%%%%%%%%%%%%%%%%%%%%%%%%%%%
\zerarcounters
\section{First order computations}
\label{sec:7}
%%%%%%%%%%%%%%%%%%%%%%%%%%%%%%%%%%%%%%%%%%%%%%%%%%%%%%%%%%%%%%%%%%%%%%%%%
%%%%%%%%%%%%%%%%%%%%%%%%%%%%%%%%%%%%%%%%%%%%%%%%%%%%%%%%%%%%%%%%%%%%%%%%%

Let us also expand the initial conditions in $\mu$:
\begin{equation}
u(0) := \bar u = u_{0}(0) + \sum_{k=1}^{\io} \mu^{k} \bar u_{k} ,
\qquad \dot u(0) := \bar v = \sum_{k=1}^{\io} \mu^{k} \bar v_{k} ,
\label{eq:7.1}
\end{equation}
and set $\Psi_{1}(\tau)=H_{1}(u_{0}(\tau),v_{0}(\tau),
\tau+\tau_{0})$ --- cf.
(\ref{eq:5.7}). We look for a solution $(u(\tau),v(\tau))$ which
is analytic in $\mu$, i.e. $u(\tau)= u_{0}(\tau) + \mu u_{1}(\tau) +
\mu^{2} u_{2}(\tau) + \ldots$ and $v(\tau)=\dot u(\tau)$ ---
cf. \cite{BBDGG,BDG,GBD3} for similar situations.
Here we are interested in the dynamics on the attractor,
hence in periodic solutions, but in principle we could also study
the dynamics near the attractor, by looking for solutions of the
form $u(\tau)=U({\rm e}^{-f_{1} \tau},{\rm e}^{-f_{2} \tau},\tau)$,
as in \cite{G1,G2,BGG}, with $f_{1}=f_{0}+O(\mu)$ and $f_{2}=O(\mu)$,
and $U(\cdot,\cdot,\psi)$ $2\pi/p$-periodic in $\psi$.
To first order one has
\begin{equation}
\left( \begin{matrix} u_{1}(\tau) \\ v_{1}(\tau) \end{matrix} \right) =
W(\tau) \left[ \left( \begin{matrix} \bar u_{1} \\ \bar v_{1} \end{matrix}
\right) + \int_{0}^{\tau} {\rm d}\tau' \, W^{-1}(\tau') \left(
\begin{matrix} 0 \\ \Psi_{1}(\tau') \end{matrix} \right) \right] ,
\label{eq:7.2}
\end{equation}
and we can confine ourselves to the first component
$u_{1}(\tau)$, since $v_{1}(\tau)=\dot u_{1}(\tau)$, 
\begin{equation}
u_{1}(\tau) = w_{11}(\tau) \,\bar u_{1} + w_{12}(\tau) \,\bar v_{1} +
\int_{0}^{\tau} {\rm d}\tau' \, {\rm e}^{F(\tau')}
\left[ w_{12}(\tau)w_{11}(\tau') - w_{11}(\tau)w_{12}(\tau') \right]
\Psi_{1}(\tau) . \nonumber
\end{equation}
which can be more conveniently written as
\begin{equation}
u_{1}(\tau) = w_{11}(\tau) \left( \bar u_{1} \!\! - \int_{0}^{\tau}
\!\!\!\! {\rm d}\tau' \, {\rm e}^{F(\tau')}
w_{12}(\tau') \Psi_{1}(\tau') \right) +
w_{12}(\tau) \left( \bar v_{1} + \!\!\int_{0}^{\tau} \!\!\!\!{\rm d}\tau' \,
{\rm e}^{F(\tau')} w_{11}(\tau') \Psi_{1}(\tau') \right) . \nonumber
\end{equation}
The function ${\rm e}^{F(\tau)}w_{11}(\tau)\Psi_{1}(\tau)$ is periodic,
while ${\rm e}^{F(\tau)}w_{12}(\tau)\Psi_{1}(\tau)$ is given by
${\rm e}^{f_{0}\tau}$ times a periodic function. Therefore
we can write $w_{11}(\tau)$ and $w_{12}(\tau)$
according to (\ref{eq:6.7}), and set --- cf. Lemma \ref{lem:6.2} ---
\begin{eqnarray}
\int_{0}^{\tau} {\rm d}\tau' \, {\rm e}^{F(\tau')}
a(\tau') \Psi_{1}(\tau')
& = & {\rm e}^{f_{0}\tau} \QQ_{1}(\tau) - \QQ_{1}(0) ,
\label{eq:7.5} \\
\int_{0}^{\tau} {\rm d}\tau' \, {\rm e}^{F(\tau')}
{\rm e}^{-f_{0}\tau'} b(\tau') \Psi_{1}(\tau')
& = & \tau \QQ_{0} + \QQ_{2}(\tau) - \QQ_{2}(0),
\label{eq:7.6}
\end{eqnarray}
for some periodic functions $\QQ_{1}(\tau)$ and $\QQ_{2}(\tau)$,
and with
\begin{equation}
\QQ_{0} = \langle {\rm e}^{\tilde F} b \, \Psi_{1} \rangle .
\label{eq:7.7}
\end{equation}

Assume that we can choose the parameters in such a way that
$\QQ_{0}=0$. Then we obtain
\begin{eqnarray}
u_{1}(\tau) \!\!\! & = & \!\!\! a(\tau) \left( \bar u_{1} + c \QQ_{1}(0) -
{\rm e}^{f_{0}\tau} c \, \QQ_{1}(\tau) + c \, \bar v_{1} +
{\rm e}^{f_{0}\tau} c \, \QQ_{1}(\tau) - c \, \QQ_{1}(0) \right.
\nonumber \\
& \!\!\! + & \!\!\! \QQ_{2}(\tau) - \QQ_{2}(0) \left) \; + \;\;
{\rm e}^{-f_{0}\tau} b(\tau) \left( \bar u_{1} + c \, \QQ_{1}(0) -
{\rm e}^{f_{0}\tau} c \, \QQ_{1}(\tau) \right) \right. ,
\label{eq:7.8}
\end{eqnarray}
and if we want that (\ref{eq:7.8}) describe a periodic function,
the constant $\bar v_{1}$ can assume any value, but we need
\begin{equation}
\bar u_{1} = - c \, \QQ_{1}(0) .
\label{eq:7.9}
\end{equation}
so that (\ref{eq:7.8}) becomes
\begin{equation}
u_{1}(\tau) = a(\tau) \left( c \, \bar v_{1} - c \, \QQ_{1}(0) +
\QQ_{2}(\tau) - \QQ_{2}(0) \right) + c \, b(\tau) \QQ_{1}(\tau) ,
\label{eq:7.10}
\end{equation}
where we have used that the function ${\rm e}^{f_{0}\tau}
c\,\QQ_{1}(\tau)$ appears twice but with opposite
sign in (\ref{eq:7.8}).

%%%%%%%%%%%%%%%%%%%%%%%%%%%%%%%%%%%%%%%%%%%%%%%%%%%%%%%%%%%%%%%%%%%%%%%%%
% REMARK 3
%%%%%%%%%%%%%%%%%%%%%%%%%%%%%%%%%%%%%%%%%%%%%%%%%%%%%%%%%%%%%%%%%%%%%%%%%
\begin{rmk} \label{rmk:3}
The constant $\bar v_{1}$ is left undetermined, and we can fix it
arbitrarily, say $\bar v_{1}=0$, as we still have at our disposal
the free parameter $\tau_{0}$; cf. \cite{GBD3}, Section 2,
for an analogous discussion.
\end{rmk}
%%%%%%%%%%%%%%%%%%%%%%%%%%%%%%%%%%%%%%%%%%%%%%%%%%%%%%%%%%%%%%%%%%%%%%%%%

Therefore we can conclude that if $\QQ_{0}=0$ then we can choose
$\bar u_{1}$ according to (\ref{eq:7.9})
in such a way that up to first order there exists
a periodic solution $u_{0}(\tau) + \mu u_{1}(\tau) + O(\mu^{2})$.
In the next Section we study in detail the condition $\QQ_{0}=0$.

%%%%%%%%%%%%%%%%%%%%%%%%%%%%%%%%%%%%%%%%%%%%%%%%%%%%%%%%%%%%%%%%%%%%%%%%%
%%%%%%%%%%%%%%%%%%%%%%%%%%%%%%%%%%%%%%%%%%%%%%%%%%%%%%%%%%%%%%%%%%%%%%%%%
\zerarcounters
\section{Compatibility to first order}
\label{sec:8}
%%%%%%%%%%%%%%%%%%%%%%%%%%%%%%%%%%%%%%%%%%%%%%%%%%%%%%%%%%%%%%%%%%%%%%%%%
%%%%%%%%%%%%%%%%%%%%%%%%%%%%%%%%%%%%%%%%%%%%%%%%%%%%%%%%%%%%%%%%%%%%%%%%%

Consider the equation $\QQ_{0}=0$, which can be written as
\begin{equation}
\eps_{1} A + B_{1}(\tau_{0}) + B_{2}(\tau_{0}) +
B_{3}(\tau_{0}) = 0 ,
\label{eq:8.1}
\end{equation}
where we have defined
\begin{eqnarray}
A \!\! & := & \!\! \left\langle {\rm e}^{\tilde F}b
\left[ \left( 1 - \beta + 3 \beta u_{0}^{2} \right) \dot u_{0} +
\frac{2}{\rho\Omega_{0}} \left( \alpha u_{0} - \beta u_{0} +
\beta\,u_{0}^{3} \right) \right] \right\rangle \nonumber \\
\!\! & = & \!\! \frac{1}{2\pi p} \int_{0}^{2\pi p}
{\rm d} \tau \, {\rm e}^{\tilde F(\tau)} b(\tau)
\left[ \dot u_{0}(\tau) \, h(u_{0}(\tau)) +
\frac{2}{\rho\Omega_{0}} k ( u_{0}(\tau)) \right]
\nonumber \\
\!\! & = & \!\! \frac{1}{2\pi\rho} \int_{0}^{2\pi\rho}
{\rm d} \tau \, {\rm e}^{\tilde F(\tau)} b(\tau)
\left[ \dot u_{0}(\tau) \, h(u_{0}(\tau)) +
\frac{2}{\rho\Omega_{0}} k ( u_{0}(\tau)) \right] ,
\label{eq:8.3}
\end{eqnarray}
and
\begin{eqnarray}
B_{1}(\tau_{0}) & := & \frac{1}{2\pi p} \int_{0}^{2\pi p}
{\rm d} \tau \, {\rm e}^{\tilde F(\tau)} b(\tau) \left[
\frac{1}{\rho\Omega_{0}}
\dot u_{0}(\tau) \left( 3 u_{0}^{2}(\tau) - 1 \right) \right]
\sin (\tau+\tau_{0}) ,
\nonumber \\
B_{2}(\tau_{0}) & := & \frac{1}{2\pi p} \int_{0}^{2\pi p}
{\rm d} \tau \, {\rm e}^{\tilde F(\tau)} b(\tau)
\left[ \frac{1}{\rho^{2}\Omega_{0}^{2}} u_{0}(\tau)
\left( u_{0}^{2}(\tau) - 1 \right) \right]
\sin (\tau+\tau_{0}) ,
\nonumber \\
B_{3}(\tau_{0}) & := & \frac{1}{2\pi p} \int_{0}^{2\pi p}
{\rm d} \tau \, {\rm e}^{\tilde F(\tau)} b(\tau)
\left[ \frac{1}{\rho\Omega_{0}} u_{0}(\tau)
\left( u_{0}^{2}(\tau) - 1 \right) \right]
\cos (\tau+\tau_{0}) .
\label{eq:8.4}
\end{eqnarray}
%

%%%%%%%%%%%%%%%%%%%%%%%%%%%%%%%%%%%%%%%%%%%%%%%%%%%%%%%%%%%%%%%%%%%%%%%%%
% REMARK 4
%%%%%%%%%%%%%%%%%%%%%%%%%%%%%%%%%%%%%%%%%%%%%%%%%%%%%%%%%%%%%%%%%%%%%%%%%
\begin{rmk} \label{rmk:4}
Note that we can write (\ref{eq:8.3}) as
$A = \frac{1}{2\pi\rho} \int_{0}^{2\pi}
{\rm d} s \, Q(s)
\left[ \frac{{\rm d}U(s)}{{\rm d}s} \, h(U(s)) +
\frac{2}{\Omega_{0}} k (U(s)) \right]$,
where $U(s)$ and $Q(s)$ are $2\pi$-periodic functions, with
$U(s)=u_{0}(\rho s)$ and $Q(s)={\rm e}^{\tilde F(\rho s)}\,b(\rho s)$.
The function $U(s)$ is the $2\pi$-periodic solution
of the differential equation
${\rm d}^{2}U/{\rm d} s^{2} + \Omega_{0}^{-1} h(U)\,
{\rm d}U/{\rm d} s + \Omega_{0}^{-2} k(U) = 0$,
and $\tilde F(\rho s) = - \frac{f_{0} \, s}{\Omega_{0}} +
\frac{1}{\Omega_{0}}
\int_{0}^{s} {\rm d}s' \, h(U(s'))$,
so that the constant $A$ is of the form $A=\bar A/\rho$, with
$\bar A$ independent of $\rho$. Hence if $A\neq0$ for some
$\rho\in\QQQ$ then it is non-zero for all rational $\rho \neq 0$.
\end{rmk}
%%%%%%%%%%%%%%%%%%%%%%%%%%%%%%%%%%%%%%%%%%%%%%%%%%%%%%%%%%%%%%%%%%%%%%%%%

By expanding
$\sin (\tau + \tau_{0}) = \sin \tau \cos \tau_{0} +
\cos \tau \sin \tau_{0}$ and
$\cos (\tau + \tau_{0}) = \cos \tau \cos \tau_{0} -
\sin \tau \sin \tau_{0}$, we can rewrite (\ref{eq:8.4}) as
$B_{i}(\tau_{0}) =  B_{i1} \cos \tau_{0} + 
B_{i2} \sin \tau_{0}$ for $i=1,2,3$,
where we have introduced the constants
\begin{eqnarray}
B_{11} \!\!\! & := & \!\!\! \frac{1}{2\pi p} \int_{0}^{2\pi p}
{\rm d} \tau \, {\rm e}^{\tilde F(\tau)} b(\tau)
\left[ \frac{1}{\rho\Omega_{0}} 
\dot u_{0}(\tau) \left( 3 u_{0}^{2}(\tau) - 1 \right) \right]
\sin \tau , \nonumber \\
B_{12} \!\!\! & := & \!\!\! \frac{1}{2\pi p} \int_{0}^{2\pi p}
{\rm d} \tau \, {\rm e}^{\tilde F(\tau)} b(\tau)
\left[ \frac{1}{\rho\Omega_{0}} 
\dot u_{0}(\tau) \left( 3 u_{0}^{2}(\tau) - 1 \right) \right]
\cos \tau , \nonumber \\
B_{21} \!\!\! & := & \!\!\! \frac{1}{2\pi p} \int_{0}^{2\pi p}
{\rm d} \tau \, {\rm e}^{\tilde F(\tau)} b(\tau)
\left[ \frac{1}{\rho^{2}\Omega_{0}^{2}} u_{0}(\tau)
\left( u_{0}^{2}(\tau) - 1 \right) \right]
\sin \tau = - \rho\Omega_{0}B_{32} , \nonumber \\
B_{22} \!\!\! & := & \!\!\! \frac{1}{2\pi p} \int_{0}^{2\pi p}
{\rm d} \tau \, {\rm e}^{\tilde F(\tau)} b(\tau)
\left[ \frac{1}{\rho^{2}\Omega_{0}^{2}} u_{0}(\tau)
\left( u_{0}^{2}(\tau) - 1 \right) \right]
\cos \tau = \rho\Omega_{0}B_{31} .
\label{eq:8.9}
\end{eqnarray}
By setting $D_{1} = - \left( B_{11} + B_{21} + B_{31} \right)$ and
$D_{2} = - \left( B_{12} + B_{22} + B_{32} \right)$,
(\ref{eq:8.1}) then becomes
\begin{equation}
\eps_{1} A = \gotD_{1}(\tau_{0}) :=
D_{1} \cos \tau_{0} + D_{2} \sin \tau_{0} .
\label{eq:8.11}
\end{equation}
All constants $B_{ij}$ in (\ref{eq:8.9}) are given by the average
of a suitable function which can be written as the product
of a $2\pi/\rho\,$-periodic function times a cosine or sine function.
Consider explicitly the constant $B_{11}$; the other constants
can be discussed in the same way. We write
\begin{equation}
B_{11} = \frac{1}{2\pi p} \int_{0}^{2\pi p} {\rm d} \tau \,
K(\tau) \, \sin \tau , \qquad \text{with} \qquad
K(\tau) = \sum_{\substack{\nu\in\ZZZ \\ \nu \; {\rm even}}}
{\rm e}^{i\nu \tau/\rho} K_{\nu} =
\sum_{\nu\in \ZZZ} {\rm e}^{i2\nu \tau/\rho} K_{2\nu} ,
\nonumber
\end{equation}
as follows from Lemmas \ref{lem:6.1} and \ref{lem:6.4}. 
If we write $\sin \tau = \sum_{\sigma = \pm 1}
(\sigma/2i) {\rm e}^{i\sigma \tau}$, then
\begin{equation}
B_{11} = \sum_{\substack{\nu\in\ZZZ,\,\sigma=\pm 1 \\
2\nu+\sigma\rho =0}} \frac{\sigma K_{2\nu}}{2i} .
\label{eq:8.14}
\end{equation}
The same argument applies to the other constants, so that
we can conclude that the constants $B_{ij}$ can be
different from zero only if $\rho$ is an even integer.
If we set $\rho=p/q$ this means $q=1$ and $p=2n$, $n\in \NNN$.
Hence for all rational $\rho\notin2\NNN$ the
first order compatibility equation (\ref{eq:8.11}) gives
$\eps_{1} A = 0 $, so that either $A=0$ and $\eps_{1}$
is arbitrary or $A \neq 0$ and $\eps_{1}=0$.
An explicit calculation (cf. Appendix \ref{app:B})
shows that $A \neq 0$. Therefore for all resonances $p\!:\!q$,
with $p/q \notin 2\NNN$, frequency locking, if possible at all,
can occur only for a range of frequencies of width
at most $\mu^{2}$.

The argument above does not imply that $D_{1},D_{2} \neq 0$
for $p/q\in 2\NNN$ --- in principle there could be cancellations
in the sum (\ref{eq:8.14}). For any given resonance $p\!:\!q$,
the non-vanishing of the constants $D_{1}$ and $D_{2}$ can be checked
numerically; for instance, when $\alpha=5$ and $\beta=4$, for $p/q=2$
one finds $D_{1}=0.00735$ and $D_{2}=-0.04507$ \cite{progress}.
Therefore for $\rho =2n$, $n\in\NNN$, frequency locking occurs
for a range of frequencies of width of order $\mu$ around the value $2n$.

%%%%%%%%%%%%%%%%%%%%%%%%%%%%%%%%%%%%%%%%%%%%%%%%%%%%%%%%%%%%%%%%%%%%%%%%%
%%%%%%%%%%%%%%%%%%%%%%%%%%%%%%%%%%%%%%%%%%%%%%%%%%%%%%%%%%%%%%%%%%%%%%%%%
\zerarcounters
\section{Higher order computations and convergence}
\label{sec:10}
%%%%%%%%%%%%%%%%%%%%%%%%%%%%%%%%%%%%%%%%%%%%%%%%%%%%%%%%%%%%%%%%%%%%%%%%%
%%%%%%%%%%%%%%%%%%%%%%%%%%%%%%%%%%%%%%%%%%%%%%%%%%%%%%%%%%%%%%%%%%%%%%%%%

To extend the analysis of the previous sections to any
perturbation order, we write the solution we are looking for as
\begin{equation}
u(\tau) = \sum_{k=0}^{\io} \mu^{k} u_{k}(\tau) , \qquad
v(\tau) = \sum_{k=0}^{\io} \mu^{k} v_{k}(\tau) = \dot u(\tau) ,
\label{eq:10.1}
\end{equation}
with $(u(0),v(0))$ written according to (\ref{eq:7.1}).
Thus we find for all $k\in\NNN$
\begin{equation}
\left( \begin{matrix} u_{k}(\tau) \\ v_{k}(\tau) \end{matrix} \right) =
W(\tau) \left[ \left( \begin{matrix} \bar u_{k} \\ \bar v_{k} \end{matrix}
\right) + \int_{0}^{\tau} {\rm d}\tau' \, W^{-1}(\tau') \left(
\begin{matrix} 0 \\ \Psi_{k}(\tau') \end{matrix} \right) \right] ,
\label{eq:10.2}
\end{equation}
where
\begin{equation}
\Psi_{k}(\tau) := \left[ \sum_{k'=1}^{k} \mu^{k'}
H_{k'} ( u(\tau), \dot u(\tau) , \tau + \tau_{0}) \right]_{k} ,
\label{eq:10.3}
\end{equation}
with $H_{k}$ defined in (\ref{eq:5.6}). The notation $[\cdot]_{k}$
for $\Psi_{k}(\tau)$ in (\ref{eq:10.3}) means the following.
In each term $H_{k'}$, we expand $u(\tau)$ and $\dot u(\tau)$
according to (\ref{eq:10.1}), and, by taking the Taylor series
of the function $H_{k'}$, we keep all contributions
proportional to $\mu^{k}$: we write the sum of these contributions
as $\mu^{k} \Psi_{k}(\tau)$. For instance one has
$\Psi_{2}(\tau) =
H_{2}(u_{0}(\tau),v_{0}(\tau),\tau + \tau_{0}) +
\frac{\partial}{\partial u_{0}}
H_{1}(u_{0}(\tau),v_{0}(\tau),\tau + \tau_{0}) \, u_{1}(\tau) +
\frac{\partial}{\partial v_{0}}
H_{1}(u_{0}(\tau),v_{0}(\tau),\tau + \tau_{0}) \, \dot u_{1}(\tau)$, 
with $u_{1}(\tau)$ given by (\ref{eq:7.10}).

As in Section \ref{sec:7}, we study only the equation for
the first component, which is
\begin{equation}
u_{k}(\tau) = w_{11}(\tau) \,\bar u_{k} + w_{12}(\tau) \,\bar v_{k} +
\int_{0}^{\tau} {\rm d}\tau' \, {\rm e}^{F(\tau')}
\left[ w_{12}(\tau)w_{11}(\tau') - w_{11}(\tau)w_{12}(\tau') \right]
\Psi_{k}(\tau') .
\label{eq:10.5}
\end{equation}
The equation (\ref{eq:10.5}) for $k=1$ has been studied in
Section \ref{sec:7}. Here we want to show that the equation
(\ref{eq:10.5}) is well defined to any perturbation order $k$,
and that it is possible to choose the constant $\eps_{k}$
in (\ref{eq:5.5}) so that it admits a periodic solution $u_{k}(\tau)$.

The discussion proceeds as in Section \ref{sec:7}, once we note that
each function $H_{k}(u,\dot u,\tau + \tau_{0})$ in (\ref{eq:5.6}) contains
a term $\eps_{k} \left( 1 - \beta + 3 \beta u^{2} \right) \dot u +
2\eps_{k}(\rho\Omega_{0})^{-1} \left[ \left( \alpha - \beta \right)
u + \beta\,u^{3} \right]$,
whereas all the other terms depend on the constants $\eps_{k'}$,
with $k'$ strictly less than $k$. Therefore for $k\in \NNN$ one has 
\begin{equation}
\Psi_{k}(\tau) = 
\eps_{k} \left( 1 - \beta + 3 \beta u_{0}^{2} \right) \dot u_{0}(\tau) +
\frac{2\eps_{k} }{\rho\Omega_{0}} \left[ \left( \alpha - \beta \right)
u_{0}(\tau) + \beta\,u_{0}^{3}(\tau) \right] + \Xi_{k}(\tau;\tau_{0}) ,
\label{eq:10.7}
\end{equation}
with the function $\Xi_{k}(\tau;\tau_{0})$ depending only on the
constants $\eps_{1},\ldots,\eps_{k-1}$, besides
the parameter $\tau_{0}$ and time $\tau$.
Therefore to any perturbation order $k$, in order to have a
periodic solution, we need
\begin{equation}
\QQ_{k,0} := \langle {\rm e}^{\tilde F} b \, \Psi_{k} \rangle = 0 ,
\label{eq:10.8}
\end{equation}
and this can be obtained by requiring
\begin{equation}
\eps_{k} A = \gotD_{k}(\tau_{0}) , \qquad
\gotD_{k}(\tau_{0}) := - \langle {\rm e}^{\tilde F} b
\, \Xi_{k}(\cdot;\tau_{0}) \rangle ,
\label{eq:10.9}
\end{equation}
with $A$ defined as in (\ref{eq:8.3}). Since $A\neq0$ (as proved in
Appendix \ref{app:B}) then we can use (\ref{eq:10.9}) to fix
$\eps_{k}$ as a function of $\tau_{0}$.
Defining the periodic functions $\QQ_{k,1}(\tau)$
and $\QQ_{k,2}(\tau)$ such that
\begin{eqnarray}
\int_{0}^{\tau} {\rm d}\tau' \, {\rm e}^{F(\tau')}
a(\tau') \Psi_{k}(\tau')
& = & {\rm e}^{f_{0}\tau} \QQ_{k,1}(\tau) - \QQ_{k,1}(0) ,
\label{eq:10.10} \\
\int_{0}^{\tau} {\rm d}\tau' \, {\rm e}^{F(\tau')}
{\rm e}^{-f_{0}\tau'} b(\tau') \Psi_{k}(\tau')
& = & \tau \QQ_{k,0} + \QQ_{k,2}(\tau) - \QQ_{k,2}(0),
\label{eq:10.11}
\end{eqnarray}
choosing the constants $\bar u_{k}$ so that
$\bar u_{k} + c \QQ_{k,1}(0) = 0$,
and using (\ref{eq:10.8}), then (\ref{eq:10.5}) gives
\begin{equation}
u_{k}(\tau) = a(\tau) \left( c \, \bar v_{k} - c \, \QQ_{k,1}(0) +
\QQ_{k,2}(\tau) - \QQ_{k,2}(0) \right) + c \, b(\tau) \QQ_{k,1}(\tau) ,
\label{eq:10.12}
\end{equation}
with the constants $\bar v_{k}$ which will be fixed
in the most convenient way (cf. Remark \ref{rmk:3}). For instance
we can set $\bar v_{k}=0$ for all $k\in\NNN$.

We can make the perturbative analysis of the previous sections rigorous
to all orders, by following the strategy introduced in \cite{GBD3,BDG},
and hence study the convergence of the perturbation series.
Alternatively, one could try to apply arguments based on
the implicit function theorems. Typically, the latter
would allow a simplification of the proof of existence
of the periodic solutions, but would be less suitable for
explicitly constructing the solutions themselves within
any given accuracy; see the comments in \cite{GBD3};
therefore we follow the first method. Note that
we are not confining ourselves to approximate
analytical solutions, which could be unreliable because
of the uncontrolled truncation of the series expansion.
On the contrary we want also to settle the issue of convergence.
In some sense this approach is complementary to that of \cite{GH},
where qualitative geometric methods are preferred to
quantitative analytical ones.

The study of the convergence of the series is standard, and
it has been discussed extensively and in full detail
in \cite{GBD3} for a similar situation. Thus, we only
sketch how the argument proceeds.

By expanding the functions $u(\tau)$ and $\dot u(\tau)$ in
$H_{k'}(u(\tau),\dot u(\tau),\tau+\tau_{0})$ in (\ref{eq:10.3})
according to (\ref{eq:10.1}), one sees that $\Psi_{k}(\tau)$
can be expressed in terms of the functions $u_{k'}(\tau)$ with $k'<k$.
On the other hand, by (\ref{eq:10.12}), the functions $u_{k}(\tau)$
are expressed in terms of the functions $\QQ_{k,1}(\tau)$
and $\QQ_{k,2}(\tau)$, which in turn
are integrals of functions involving $\Psi_{k}(\tau)$,
and hence depend on $u_{k'}(\tau)$ for $k'<k$.

This means that we have recursive equations for the functions
$u_{k}(\tau)$. By passing to the Fourier space, that is by expanding
$u_{k}(\tau) = \sum_{\nu\in\ZZZ} {\rm e}^{i\nu \tau/p} u_{k,\nu}$,
we obtain recursive equations for the Fourier coefficients $u_{k,\nu}$.
We do not write them explicitly because
the ensuing expressions are rather cumbersome, but one can
easily work out the analytical expressions for the recursions
by following the scheme that we have outlined.
Eventually, we can represent $u_{k,\nu}$
for $k\ge 1$ and $\nu\in\ZZZ$, in terms of tree graphs,
which can be studied with the techniques of \cite{GBD3}.

We do not repeat the analysis here, but we instead just
give the final result. To any order $k\ge 1$ 
one obtain the following bounds for the Fourier coefficients:
$| u_{k,\nu}| \le C_{1} C_{2}^{k-1}$ and
$\sum_{\nu\in\ZZZ}| u_{k,\nu}| \le C_{3} C_{2}^{k-1}$,
for suitable positive constants $C_{1},C_{2},C_{3}$,
depending on $\rho$. This implies the convergence of
the perturbation series (\ref{eq:10.1}) for $\mu$ small enough,
say for $|\mu|<C_{2}^{-1}$.

%%%%%%%%%%%%%%%%%%%%%%%%%%%%%%%%%%%%%%%%%%%%%%%%%%%%%%%%%%%%%%%%%%%%%%%%%
%%%%%%%%%%%%%%%%%%%%%%%%%%%%%%%%%%%%%%%%%%%%%%%%%%%%%%%%%%%%%%%%%%%%%%%%%
\zerarcounters
\section{Arnold tongues and devil's staircase}
\label{sec:11}
%%%%%%%%%%%%%%%%%%%%%%%%%%%%%%%%%%%%%%%%%%%%%%%%%%%%%%%%%%%%%%%%%%%%%%%%%
%%%%%%%%%%%%%%%%%%%%%%%%%%%%%%%%%%%%%%%%%%%%%%%%%%%%%%%%%%%%%%%%%%%%%%%%%

We use the perturbative analysis, developed to all orders in the
previous section, to study for which values of the driving frequency
$\omega$ one has locking. We shall see that the analysis accounts
for the devil's staircase structure found in \cite{OBYK},
for small values of the driving amplitude $\mu$.

%%%%%%%%%%%%%%%%%%%%%%%%%%%%%%%%%%%%%%%%%%%%%%%%%%%%%%%%%%%%%%%%%%%%%%%%%
% LEMMA 7
%%%%%%%%%%%%%%%%%%%%%%%%%%%%%%%%%%%%%%%%%%%%%%%%%%%%%%%%%%%%%%%%%%%%%%%%%
\begin{lemma} \label{lem:11.1}
The functions $H_{k}(u,\dot u,\tau + \tau_{0})$ in (\ref{eq:5.6})
are polynomials of odd order in $(u,\dot u)$ for all $k\in\NNN$.
\end{lemma}
%%%%%%%%%%%%%%%%%%%%%%%%%%%%%%%%%%%%%%%%%%%%%%%%%%%%%%%%%%%%%%%%%%%%%%%%%

%%%%%%%%%%%%%%%%%%%%%%%%%%%%%%%%%%%%%%%%%%%%%%%%%%%%%%%%%%%%%%%%%%%%%%%%%
\prova The function $H(u,\dot u,\ddot u,\mu)$ given by (\ref{eq:4.1})
is a polynomial of odd order in $(u,\dot u,\ddot u)$. By writing
$H(u,\dot u,\ddot u,\mu)$ as in (\ref{eq:5.6}), the only
term containing $\ddot u$ is the first one ($k=0$), so that
all the other terms are polynomials of odd order in $(u,\dot u)$. \EP
%%%%%%%%%%%%%%%%%%%%%%%%%%%%%%%%%%%%%%%%%%%%%%%%%%%%%%%%%%%%%%%%%%%%%%%%%

%%%%%%%%%%%%%%%%%%%%%%%%%%%%%%%%%%%%%%%%%%%%%%%%%%%%%%%%%%%%%%%%%%%%%%%%%
% LEMMA 8
%%%%%%%%%%%%%%%%%%%%%%%%%%%%%%%%%%%%%%%%%%%%%%%%%%%%%%%%%%%%%%%%%%%%%%%%%
\begin{lemma} \label{lem:11.2}
For all $k\in\NNN$ one has
\begin{eqnarray}
u_{k}(\tau) &\!\!\! = \!\!\! &
\sum_{\substack{\nu\in\ZZZ\\\nu \; {\rm odd}}}
\sum_{\substack{\s \in \ZZZ \\ |\s| \le k}}
{\rm e}^{i\nu \tau/\rho} {\rm e}^{i\s (\tau + \tau_{0})}
\overline u_{k,\nu,\sigma}
\label{eq:11.1} \\
\Psi_{k}(\tau) &\!\!\! = \!\!\! &
\sum_{\substack{\nu\in\ZZZ\\\nu \; {\rm odd}}}
\sum_{\substack{\s \in \ZZZ \\ |\s| \le k}}
{\rm e}^{i\nu \tau/\rho} {\rm e}^{i\s (\tau + \tau_{0})} \overline
\Psi_{k,\nu,\sigma} ,
\label{eq:11.2} \end{eqnarray}
with the coefficients $\overline u_{k,\nu,\sigma}$
and $\overline \Psi_{k,\nu,\sigma}$ independent of $\tau_{0}$.
\end{lemma}
%%%%%%%%%%%%%%%%%%%%%%%%%%%%%%%%%%%%%%%%%%%%%%%%%%%%%%%%%%%%%%%%%%%%%%%%%

%%%%%%%%%%%%%%%%%%%%%%%%%%%%%%%%%%%%%%%%%%%%%%%%%%%%%%%%%%%%%%%%%%%%%%%%%
\prova First of all note that if $\Psi_{k}(\tau)$ is  of the form
(\ref{eq:11.2}) then $u_{k}(\tau)$ is also of the form (\ref{eq:11.1}).
This can be proved as follows. For brevity, here and henceforth
we say that $u_{k}(\tau)$ and $\Psi_{k}(\tau)$ `contain only odd
harmonics' if they are of the form (\ref{eq:11.1}) and
(\ref{eq:11.2}), respectively.
The functions $\QQ_{k,1}(\tau)$ and $\QQ_{k,2}(\tau)$ are integrals
of functions which are either periodic functions $P(\tau)$ or
of the form ${\rm e}^{f_{0}t}$ times periodic functions $P(\tau)$.
In all cases the function $P(\tau)$ is given by the product of three
functions: two of these functions --- one is either $a(\tau)$ or $b(\tau)$,
the other one is $\Psi_{k}(\tau)$ --- contain odd harmonics,
by Lemma \ref{lem:6.4} and by our assumption
on $\Psi_{k}(\tau)$, while the third one --- ${\rm e}^{\tilde F(\tau)}$ ---
contains only even harmonics. If we compare (\ref{eq:6.11})
with (\ref{eq:6.12}) we see that the integral of a function
${\rm e}^{C \tau} P(\tau)$ is of the form $D+{\rm e}^{C\tau }Q(\tau)$,
where $Q(\tau)$ contains the same harmonics as $P(\tau)$. Therefore
both $\QQ_{k,1}(\tau)$ and $\QQ_{k,2}(\tau)$
are periodic functions containing only even harmonics.
Then, recall that $u_{k}(\tau)$ is given by (\ref{eq:10.12}).
We have already used the fact that the functions $a(\tau)$ and $b(\tau)$
contain only odd harmonics, so that we can conclude that, as claimed
above, if $\Psi_{k}(\tau)$ is  of the form (\ref{eq:11.2}) then
$u_{k}$ is of the form (\ref{eq:11.1}).

Then, the proof of the lemma proceeds by induction. Recall that for
$k=1$ one has $\Psi_{1}(\tau)=H_{1}(u_{0}(\tau),\dot u_{0}(\tau),
\tau + \tau_{0})$, with $H_{1}$ given by (\ref{eq:5.7}), so that, by
Lemma \ref{lem:3.1} and Lemma \ref{lem:11.1}, $\Psi_{1}(\tau)$
is of the form (\ref{eq:11.1}), and, by the previous observation,
the function $u_{1}(\tau)$ is also of the form (\ref{eq:11.2}).

By assuming that $u_{k}(\tau)$ is of the form (\ref{eq:11.1}) for all
$k<\bar k$, then by Lemma \ref{lem:11.1} it also follows that
$\Psi_{\bar k}(\tau)$, given by (\ref{eq:11.2}), is of the form
(\ref{eq:11.2}). Again by the observation at the beginning
of the proof, it follows that $u_{\bar k}(\tau)$
can be expressed as in (\ref{eq:11.2}). \EP
%%%%%%%%%%%%%%%%%%%%%%%%%%%%%%%%%%%%%%%%%%%%%%%%%%%%%%%%%%%%%%%%%%%%%%%%%

%%%%%%%%%%%%%%%%%%%%%%%%%%%%%%%%%%%%%%%%%%%%%%%%%%%%%%%%%%%%%%%%%%%%%%%%%
% REMARK 5
%%%%%%%%%%%%%%%%%%%%%%%%%%%%%%%%%%%%%%%%%%%%%%%%%%%%%%%%%%%%%%%%%%%%%%%%%
\begin{rmk} \label{rmk:11.2}
If we expand $u_{k}(\tau)$ as a Fourier series,
$ u_{k}(\tau) = \sum_{\nu\in\ZZZ} {\rm e}^{i\nu \tau/p} u_{k,\nu}$,
then (\ref{eq:11.1}) implies
$$ u_{k,\nu} = \sum_{\substack{\nu'\in\ZZZ,|\s|\le k \\ q\nu'+p \s 
= \nu}} {\rm e}^{i\s \tau_{0}} \overline u_{k,\nu',\s} . $$
In particular $u_{k}(\tau)$ and $\Psi_{k}(\tau)$ are polynomials
of order $k$ in $\tau_{0}$.
\end{rmk}
%%%%%%%%%%%%%%%%%%%%%%%%%%%%%%%%%%%%%%%%%%%%%%%%%%%%%%%%%%%%%%%%%%%%%%%%%

%%%%%%%%%%%%%%%%%%%%%%%%%%%%%%%%%%%%%%%%%%%%%%%%%%%%%%%%%%%%%%%%%%%%%%%%%
% LEMMA 9
%%%%%%%%%%%%%%%%%%%%%%%%%%%%%%%%%%%%%%%%%%%%%%%%%%%%%%%%%%%%%%%%%%%%%%%%%
\begin{lemma} \label{lem:11.3}
For all $k\in\NNN$ one has
$$ \gotD_{k}(\tau_{0}) = \gotD_{k,0} + \frac{1}{2\pi p}
\sum_{\substack{\nu\in\ZZZ \\ \nu \; {\rm even}}} \,
\sum_{\substack{\s\in\ZZZ \\ 0<|\s| \le k}}
\int_{0}^{2\pi p} {\rm d} \tau \,
{\rm e}^{i\nu \tau/\rho} {\rm e}^{i\s (\tau + \tau_{0})} K_{k,\nu,\s} , $$
for suitable $\tau_{0}$-independent coefficients $K_{k,\nu,\s}$,
depending on $\eps_{1},\ldots,\eps_{k-1}$, but not on $\eps_{k}$.
\end{lemma}
%%%%%%%%%%%%%%%%%%%%%%%%%%%%%%%%%%%%%%%%%%%%%%%%%%%%%%%%%%%%%%%%%%%%%%%%%

%%%%%%%%%%%%%%%%%%%%%%%%%%%%%%%%%%%%%%%%%%%%%%%%%%%%%%%%%%%%%%%%%%%%%%%%%
\prova The functions ${\rm e}^{\tilde F(\tau)}$ and $b(\tau)$ in
(\ref{eq:10.8}) are periodic in $\tau$ with period $2\pi\rho=2\pi p/q$,
and contain only even and odd harmonics, respectively, whereas
$\Psi_{k}(\tau)$ is given by (\ref{eq:11.2}). By Lemma \ref{lem:11.2},
this yields that $\QQ_{k,0}=\QQ_{k,0}(\tau_{0})$ is of the form
$ \QQ_{k,0} = \frac{1}{2\pi p} 
\sum_{\substack{\nu\in\ZZZ \\ \nu \; {\rm even}}}
\sum_{\substack{\s\in\ZZZ \\ |\s| \le k}}
\int_{0}^{2\pi p} {\rm d} \tau \,
{\rm e}^{i\nu \tau/\rho} {\rm e}^{i\s (\tau + \tau_{0})} Q_{k,\nu,\s}$,
for suitable coefficients $Q_{k,\nu,\s}$, which are
independent of $\tau_{0}$ but depend on $\eps_{1},\ldots,\eps_{k}$.
In particular the only contribution to $\QQ_{k,0}$ depending
on $\eps_{k}$ is of the form $\eps_{k} A$ --- cf. (\ref{eq:10.9}) ---,
so that we can write $\QQ_{k,0}=\eps_{k}A + \overline \gotD_{k,0}
(\eps_{1},\ldots,\eps_{k-1};\tau_{0})$, for a suitable function
$\overline \gotD_{k,0} (\eps_{1},\ldots,\eps_{k-1};\tau_{0})$.\EP
%%%%%%%%%%%%%%%%%%%%%%%%%%%%%%%%%%%%%%%%%%%%%%%%%%%%%%%%%%%%%%%%%%%%%%%%%

By (\ref{eq:5.5}) and (\ref{eq:10.9}), and using
Lemma \ref{lem:11.3}, we can write
\begin{equation}
\eps(\mu) = \gotD(\tau_{0},\mu) :=
\frac{1}{A} \sum_{k=1}^{\io} \mu^{k} \gotD_{k}(\tau_{0}) ,
\qquad \gotD_{k}(\tau_{0}) = \sum_{\substack{\s\in\ZZZ \\
|\s| \le k}} {\rm e}^{i\s \tau_{0}} \gotD_{k,\s} ,
\label{eq:11.3} \end{equation}
for suitable coefficients $\gotD_{k,\s}$.
For given $\omega$, for a periodic solution with period $2\pi p$
to exist, we need that $\eps(\mu)$, defined according to (\ref{eq:5.1}),
satisfy (\ref{eq:11.3}) for some $\tau_{0}\in[0,2\pi)$.
Therefore, by defining
\begin{equation}
\eps_{\rm max}(\rho) :=
\max_{0 \le \tau_{0} \le 2\pi} \gotD(\tau_{0},\mu) , \qquad
\eps_{\rm min}(\rho) :=
\min_{0 \le \tau_{0} \le 2\pi} \gotD(\tau_{0},\mu) , \nonumber
\end{equation}
and setting $\WW(\rho) = \eps_{\rm max}(\rho) -
\eps_{\rm min}(\rho)$, such a periodic solution exists
for all $\eps(\mu) \in \WW(\rho)$.

%%%%%%%%%%%%%%%%%%%%%%%%%%%%%%%%%%%%%%%%%%%%%%%%%%%%%%%%%%%%%%%%%%%%%%%%%
% LEMMA 10
%%%%%%%%%%%%%%%%%%%%%%%%%%%%%%%%%%%%%%%%%%%%%%%%%%%%%%%%%%%%%%%%%%%%%%%%%
\begin{lemma} \label{lem:11.4}
Fix $\rho=p/q$. One has $\gotD_{k}(\tau_{0})=\gotD_{k,0}$
for all $k<q$ if $p$ is even and for all $k<2q$ if $p$ is odd.
\end{lemma}
%%%%%%%%%%%%%%%%%%%%%%%%%%%%%%%%%%%%%%%%%%%%%%%%%%%%%%%%%%%%%%%%%%%%%%%%%

%%%%%%%%%%%%%%%%%%%%%%%%%%%%%%%%%%%%%%%%%%%%%%%%%%%%%%%%%%%%%%%%%%%%%%%%%
\prova One can write $\eps_{k}A=\gotD_{k}(\tau_{0})$,
with $\gotD_{k}(\tau_{0})$ defined in Lemma \ref{lem:11.3}.
By comparing (\ref{eq:11.3}) with the expression for
$\gotD_{k}(\tau_{0})$ in Lemma \ref{lem:11.3}, we see that 
$\gotD_{k,\s} = \frac{1}{2\pi p}
\sum_{\substack{\nu\in\ZZZ \\ \nu \; {\rm even}}}
\int_{0}^{2\pi p} {\rm d} \tau \, {\rm e}^{i\nu \tau/\rho}
{\rm e}^{i\s (\tau + \tau_{0})} K_{k,\nu,\s}$
for $\s\neq0$, so that one can have $\gotD_{k,\s} \neq 0$ only if $\s\rho \in
2\NNN$ for some $|\s|\le k$. Hence, if $\rho=p/q$ with either
even $p$ and $q>k$ or odd $p$ and $q>2k$, one has $\gotD_{k,\s}=0$.
In other words, for fixed $\rho=p/q$ one has
$\gotD_{k}(\tau_{0})=\gotD_{k,0}$ for all $k < q$ if $p$
is even and for all $k<2q$ if $p$ is odd.\EP
%%%%%%%%%%%%%%%%%%%%%%%%%%%%%%%%%%%%%%%%%%%%%%%%%%%%%%%%%%%%%%%%%%%%%%%%%

By Lemma \ref{lem:11.4} we can write in (\ref{eq:11.3})
\begin{equation}
\gotD_{k}(\tau_{0}) = \gotD_{k,0} + \widetilde \gotD_{k}(\tau_{0}) ,
\qquad \widetilde \gotD_{k}(\tau_{0}) = \sum_{\substack{\s\in\ZZZ \\
0<|\s| \le k}} {\rm e}^{i\s \tau_{0}} \gotD_{k,\s} , \nonumber
\end{equation}
where the zero-mean function $\widetilde \gotD_{k}(\tau_{0})$
vanishes for $k<q$ if $p$ is even and
for $k<2q$ if $p$ is odd.

%%%%%%%%%%%%%%%%%%%%%%%%%%%%%%%%%%%%%%%%%%%%%%%%%%%%%%%%%%%%%%%%%%%%%%%%%
% REMARK 6
%%%%%%%%%%%%%%%%%%%%%%%%%%%%%%%%%%%%%%%%%%%%%%%%%%%%%%%%%%%%%%%%%%%%%%%%%
\begin{rmk} \label{rmk:6}
The coefficient $\gotD_{k,0}$ does not contribute to $\WW(\rho)$:
when making the difference
between $\eps_{\rm max}$ and $\eps_{\rm min}$ only
$\widetilde \gotD_{k}(\tau_{0})$ plays a role.
Therefore Lemma \ref{lem:11.4} implies that
$\WW(\rho)=O(\mu)$ only for $\rho=2n$, $n\in\NNN$,
$\WW(\rho)=O(\mu^{2})$ only for $\rho = 2n-1$, $n\in\NNN$,
$\WW(\rho)=O(\mu^{3})$ only for $\rho = 2n/3$, $n\in\NNN$,
$\WW(\rho)=O(\mu^{4})$ only for $\rho = (2n-1)/2$, $n\in\NNN$,
$\WW(\rho)=O(\mu^{5})$ only for $\rho = 2n/5$, $n\in\NNN$,
$\WW(\rho)=O(\mu^{6})$ only for $\rho = (2n-1)/3$, $n\in\NNN$,
and so on. In general, if $\rho=p/q$ with even $p$,
then $\WW(\rho)=O(\mu^{q})$, while if $\rho=p/q$ with odd $p$,
then $\WW(\rho)=O(\mu^{2q})$.
\end{rmk}
%%%%%%%%%%%%%%%%%%%%%%%%%%%%%%%%%%%%%%%%%%%%%%%%%%%%%%%%%%%%%%%%%%%%%%%%%

If we recall the definition (\ref{eq:5.1}) of $\eps(\mu)$ and we set
\begin{equation}
\omega_{\rm min}(\rho) := \frac{\rho\Omega_{0}}{1 +
\rho\Omega_{0} \, \eps_{\rm max}(\rho)} , \qquad
\omega_{\rm max}(\rho) := \frac{\rho\Omega_{0}}{1 +
\rho\Omega_{0} \, \eps_{\rm min}(\rho)} ,
\label{eq:11.6} \end{equation}
we obtain that for
\begin{equation}
\omega_{\rm min}(\rho) \le \omega \le \omega_{\rm max}(\rho)
\label{eq:11.7} \end{equation}
there exists a periodic solution with period $2\pi p$
(recall that $\rho=p/q$). In the $(\omega,\mu)$ plane 
the region (\ref{eq:11.7}) defines a distorted wedge
with apex at $\omega=\rho\Omega_{0}$ on the real axis. 

Call $\Delta\omega(\rho)=\omega_{\rm max}(\rho)-\omega_{\rm min}
(\rho)$ the range of frequencies around the value
$\rho\Omega_{0}$, with $\rho=p/q$, for which there is
frequency locking. Then
\begin{equation}
\Delta\omega(2n/k) = O(\mu^{k}) , \qquad
\Delta\omega((2n+1)/k) = O(\mu^{2k})
\label{eq:11.8} \end{equation}
for all $k,n \in \NNN$ such that $2n/k$ and $(2n+1)/k$, respectively,
are irreducible fractions. Indeed, $\Delta\omega(\rho)$ is proportional
to $\WW(\rho)$, so that $\gotD_{k,0}$ does not contribute to
the width of the plateau, but only to its `centre'.
In the $(\omega,\mu)$ plane the locking regions (Arnold tongues)
`emanate' from  the values $\rho\Omega_{0}$, with $\rho\in\QQQ$. For
$\rho\in2\NNN$ they are centred around the vertical passing through
$\omega=\rho\Omega_{0}$ and for fixed $\mu$ have width $O(\mu)$.
For all the other rational values of $\rho$, in general, they
slightly bend away from the vertical: for fixed $\mu$
the centre of the region is shifted of order $\mu^{2}$
with respect to the value $\omega=p\Omega_{0}/q$, whereas the width
is $O(\mu^{q})$ for even $p$ and $O(\mu^{2q})$ for odd $p$.

%%%%%%%%%%%%%%%%%%%%%%%%%%%%%%%%%%%%%%%%%%%%%%%%%%%%%%%%%%%%%%%%%%%%%%%%%
%%%%%%%%%%%%%%%%%%%%%%%%%%%%%%%%%%%%%%%%%%%%%%%%%%%%%%%%%%%%%%%%%%%%%%%%%
\zerarcounters
\section{Conclusions and open problems}
\label{sec:13}
%%%%%%%%%%%%%%%%%%%%%%%%%%%%%%%%%%%%%%%%%%%%%%%%%%%%%%%%%%%%%%%%%%%%%%%%%
%%%%%%%%%%%%%%%%%%%%%%%%%%%%%%%%%%%%%%%%%%%%%%%%%%%%%%%%%%%%%%%%%%%%%%%%%

The locking of oscillators onto subharmonics of the driving frequency
(also called frequency demultiplication) has been well
known in electronics since the work of van der Pol and van der Mark
\cite{PM}. In the $(\omega,\mu)$ frequency-amplitude plane,
the locking region occurs in distorted wedges (Arnold tongues)
with apices corresponding to the rational values on the frequency axis.
If one plots the ratio of the driver frequency $\omega$ to
the output frequency $\Omega$ versus the driving frequency $\omega$,
one obtains a so-called devil's staircase, i.e.
a self-similar fractal object, where the qualitative
structure is replicated at a higher level of resolution,
with plateaux corresponding to rational values of the ratio.

The phase locking phenomenon, the existence of the Arnold tongues,
and the devil's staircase picture have been proved rigorously in
some mathematical models, such as the circle map \cite{A},
and studied numerically for several electronic circuits,
such as the van der Pol equation \cite{GH}, the Josephson junction
\cite{AC,L1,QWZ}, the Chua circuit \cite{PZC} among others.

In this paper we have studied analytically the
injection-locked frequency divider equation considered in \cite{OBYK}.
In particular we aimed to understand the devil's staircase picture, 
with the largest plateaux corresponding to integer resonances of
even order, and to provide an algorithm to compute the width of
the plateaux for small values of the driving amplitude $\mu$.

The main result is summarised by (\ref{eq:11.8}),
which gives the width of the Arnold tongues in terms of the
driving amplitude $\mu$ and of the resonances $p\!:\!q$.
Note that the width of the tongues is narrower
for resonances of higher order.

In most of the analytic discussions in the
literature, one usually assumes that the unperturbed system
is written in a very simple form --- see for instance \cite{GH}.
Of course, determining analytically the change of variables
which puts the system into such a form can be very difficult
in general, in principle as difficult as finding explicitly 
the solution itself. Hence, we have
preferred to work directly with the original coordinates.
Even if we have concentrated here on the
injection-locked frequency divider equation, our analysis
applies to any driven Li\'enard equation, of which
the van der Pol equation is a particular type.
The dynamics of the forced or driven van der Pol equation has been
analytically investigated in \cite{L3,L4,L2}. However,
we could not rely on results existing in the literature,
as we are interested in the exact structure of the Arnold
tongues, which of course strongly depends on the
particular form of the system under study.

We have considered the model (\ref{eq:2.1})
introduced in \cite{OBYK}. In particular we have taken
the same driving term as in \cite{OBYK}, containing only
one non-zero harmonic. In principle, one can consider
more general functions, for instance any analytic periodic function,
instead of the sine function. In that case the driving function
contains all the harmonics; of course, by analyticity,
the coefficients of the harmonics decay exponentially fast.
Then one could ask how the analysis changes in such a case.
From a technical point of view, there are no further complications.
However, the conclusions about the devil's staircase structure
are slightly different.
For instance, the width of all plateaux becomes of order $\mu$.
This follows by the same arguments as given in Section \ref{sec:8}.
The analogues of the functions $B_{i}(\tau_{0})$ in (\ref{eq:8.4})
contain all the harmonics $\sin(\sigma(\tau+\tau_{0}))$
and $\cos(\sigma(\tau+\tau_{0}))$, with $\sigma\in\NNN$, so that,
when imposing the constraint $2\nu+\sigma\rho=0$ in (\ref{eq:8.14}),
one no longer has $\s=\pm 1$. On the contrary, one has
$\sigma\in\ZZZ$; thus in general the constraint can be satisfied
for all $\rho\in\QQQ$ (by choosing $\nu$ appropriately),
and so all the plateaux have width of order $\mu$.
However, the larger $p$ and $q$ in $\rho=p/q$ are, the narrower
the plateau is: indeed $2\nu+\s \rho=(2 \nu q + \s p)/q=0$
requires $2|\nu|/|\sigma|=p/q$, hence, for very large values
of $p$ and $q$, both $\nu$ and $\sigma$ are very large,
and hence the factors $K_{2\nu}$ contributing to $B_{11}$
in (\ref{eq:8.14}) are very small. This is consistent
with the fact that the union of Arnold tongues
form an open dense subset of the $(\omega,\mu)$ plane,
whose complement converges to full
measure as $\mu\to0$ \cite{He}.
So, an important observation is that large plateaux
have not been found in \cite{OBK} for odd integer values
because of the peculiar form of the driving term: they
would appear by taking, for instance, a
driving term involving also the harmonics with $\nu=\pm 2$.

We have studied analytically the existence
and properties of the periodic solution which continues
the unperturbed limit cycle when the perturbation is switched on.
It would be interesting to prove analytically also that such a solution
is attracting, for instance by determining the Lyapunov exponents
or studying the more general solutions which move
nearby and tend asymptotically to the attractor --- for instance
by following the strategy outlined in the first paragraph
of Section \ref{sec:7}.

Another interesting problem to investigate analytically
concerns the dynamics far away from the resonances,
i.e. when the rotation vector $(\omega,\Omega_{0})$
satisfies some Diophantine condition such as
the standard Diophantine condition mentioned
in Section \ref{sec:1} --- see also the comments
in the last paragraph of Section \ref{sec:2} --- or
the weaker Bryuno condition \cite{G3,GBD2}.
Such values of $\omega$, in the devil's staircase picture,
are complementary to those for which frequency locking occurs.

The analysis we have performed is based on perturbation theory,
and applies for $\mu$ small enough. It would be interesting
to investigate the locking diagram in the
$(\om,\mu)$ plane for large values of $\mu$.
It could be worthwhile to enquire further both analytically
(for small values of $\mu$) and numerically (even for larger
values of $\mu$) into the structure of the Arnold tongues in the
$(\omega,\mu)$ plane. Work is underway concerning these
problems \cite{progress}.

%%%%%%%%%%%%%%%%%%%%%%%%%%%%%%%%%%%%%%%%%%%%%%%%%%%%%%%%%%%%%%%%%%%%%%%%%
\vspace{.5truecm}
\noindent \textbf{Acknowledgements.} We thank Giovanni Gallavotti
for useful discussions, and Peter Kennedy for 
bringing this problem to our attention. We are also indebted to
Henk Bruin and Freddy Dumortier
for providing us with the references \cite{C} and \cite{Z}.

%%%%%%%%%%%%%%%%%%%%%%%%%%%%%%%%%%%%%%%%%%%%%%%%%%%%%%%%%%%%%%%%%%%%%%%%%
\appendix
%%%%%%%%%%%%%%%%%%%%%%%%%%%%%%%%%%%%%%%%%%%%%%%%%%%%%%%%%%%%%%%%%%%%%%%%%

%%%%%%%%%%%%%%%%%%%%%%%%%%%%%%%%%%%%%%%%%%%%%%%%%%%%%%%%%%%%%%%%%%%%%%%%%
%%%%%%%%%%%%%%%%%%%%%%%%%%%%%%%%%%%%%%%%%%%%%%%%%%%%%%%%%%%%%%%%%%%%%%%%%
\zerarcounters
\section{Well-posedness of the Wronskian matrix}
\label{app:A}
%%%%%%%%%%%%%%%%%%%%%%%%%%%%%%%%%%%%%%%%%%%%%%%%%%%%%%%%%%%%%%%%%%%%%%%%%
%%%%%%%%%%%%%%%%%%%%%%%%%%%%%%%%%%%%%%%%%%%%%%%%%%%%%%%%%%%%%%%%%%%%%%%%%

Let $u_{0}$ be the periodic solution of (\ref{eq:5.2}) satisfying the
conditions (\ref{eq:6.3}). Write $u_{0}(\tau)=r_{0}+r_{1}\tau^{2}/2+
r_{2}\tau^{3}/3+O(\tau^{4})$ --- cf. Remark \ref{rmk:1}.

%%%%%%%%%%%%%%%%%%%%%%%%%%%%%%%%%%%%%%%%%%%%%%%%%%%%%%%%%%%%%%%%%%%%%%%%%
\begin{lemma} \label{lem:A.1}
The function $w_{11}(\tau)$ in (\ref{eq:6.7}) is smooth.
\end{lemma}
%%%%%%%%%%%%%%%%%%%%%%%%%%%%%%%%%%%%%%%%%%%%%%%%%%%%%%%%%%%%%%%%%%%%%%%%%

%%%%%%%%%%%%%%%%%%%%%%%%%%%%%%%%%%%%%%%%%%%%%%%%%%%%%%%%%%%%%%%%%%%%%%%%%
\prova By deriving (\ref{eq:5.2}), one finds
\begin{equation}
\dddot u_{0} + f(u_{0})\, \ddot u_{0} + f'(u_{0})\,\dot u_{0}^{2} +
g'(u_{0})\,\dot u_{0} = 0 ,
\label{eq:A1}
\end{equation}
where $f'$ and $g'$ are the derivatives of $f$ and $g$
with respect to their arguments, while the dots denote derivatives
with respect to the time $\tau$.

By computing (\ref{eq:A1}) at $\tau=0$ and using
that $\dot u_{0}(0)=0$, we find
\begin{equation}
0 = \dddot u_{0}(0) + f(u_{0}(0))\, \ddot u_{0}(0) 
= 2r_{2} + f(r_{0})\,r_{1} .
\label{eq:A2}
\end{equation}
In (\ref{eq:6.8}) we can write $F(\tau) = \int_{0}^{\tau}
{\rm d}\tau' \, f(u_{0}(0)) + O(\tau^{2}) =
f(r_{0}) \, \tau + O(\tau^{2})$,
so that ${\rm e}^{-F(\tau)} = 1 - f(r_{0}) \, \tau + O(\tau^{2})$.
On the other hand one has
$1/\dot u_{0}^{2}(\tau) = (r_{1}^{2} \tau^{2})^{-1} \left( 1 -
2r_{2}\tau/r_{1} + O(\tau^{2}) \right)$.
Therefore the integrand in (\ref{eq:6.7}) can be expanded as
\begin{equation}
\frac{{\rm e}^{-F(\tau)}}{\dot u_{0}^{2}(\tau)} =
\frac{1}{r_{1}^{2} \tau^{2}} \left( 1 - \frac{2r_{2}}{r_{1}} \tau -
f(r_{0}) \, \tau + O(\tau^{2}) \right) .
\label{eq:A6}
\end{equation}
The term $1/r_{1}^{2} \tau^{2}$ produces a linear divergence,
which is compensated by the function $\dot u_{0}(\tau)$ in front of the
integral. The integral arising from the linear term inside the
parentheses of (\ref{eq:A6}) would produce a logarithmic divergence
(hence a divergence of the first derivative of $w_{11}(\tau)$); however
such a term is of the form $- \left( 2r_{2}/r_{1} + f(r_{0}) \right) \tau =
- \tau \, r_{1}^{-1} \left( 2r_{2} + f(r_{0}) \, r_{1} \right)$,
which vanishes because of (\ref{eq:A2}). Finally, the remaining
part of the integrand arises from the terms
of order $\tau^{2}$ in (\ref{eq:A6}), and hence produces regular terms.
This proves that the function $w_{11}(\tau)$ is smooth.\EP
%%%%%%%%%%%%%%%%%%%%%%%%%%%%%%%%%%%%%%%%%%%%%%%%%%%%%%%%%%%%%%%%%%%%%%%%%

%%%%%%%%%%%%%%%%%%%%%%%%%%%%%%%%%%%%%%%%%%%%%%%%%%%%%%%%%%%%%%%%%%%%%%%%%
\begin{lemma} \label{lem:A.2}
There exists a unique $\bar\tau \in (0,\pi\rho)$
such that $\dot w_{11}(0)=0$.
\end{lemma}
%%%%%%%%%%%%%%%%%%%%%%%%%%%%%%%%%%%%%%%%%%%%%%%%%%%%%%%%%%%%%%%%%%%%%%%%%

%%%%%%%%%%%%%%%%%%%%%%%%%%%%%%%%%%%%%%%%%%%%%%%%%%%%%%%%%%%%%%%%%%%%%%%%%
\prova One can write $w_{11}(\tau)$ in (\ref{eq:6.7}) as $w_{11}(\tau)
= c_{1} \dot u_{0}(\tau) \left( R(\tau) - R(\bar\tau) \right)$,
where $R(\tau)$ is a primitive of the function
${\rm e}^{-F(\tau)}/\dot u_{0}^{2}(\tau)$, i.e.
$\dot R(\tau) = r(\tau) := {\rm e}^{-F(\tau)}/\dot u_{0}^{2}(\tau)$.

The function $r(\tau)$ is smooth and strictly positive
for $t\in (0,\pi\rho)$, and hence its primitive $R(\tau)$ is strictly
increasing for $t\in (0,\pi\rho)$. For all $t,\bar\tau \in (0,\pi\rho)$
the function
\begin{equation}
R(\tau,\bar\tau) := \int_{\bar\tau}^{\tau} {\rm d}\tau' \, r(\tau') =
R(\tau) - R(\bar\tau)
\label{eq:A10}
\end{equation}
is smooth, and for all $\bar\tau \in (0,\pi\rho)$ one has
$\lim_{\tau\to 0^{+}} R(\tau,\bar\tau) = - \infty$ and
$\lim_{\tau\to \pi\rho^{-}} R(\tau,\bar\tau) = + \infty$,
which imply that for all $\bar\tau \in (0,\pi\rho)$ the
function $R(\tau,\bar\tau)$ is strictly increasing in $\tau$
from $-\infty$ to $+\infty$. Now $\dot w_{11}(\tau) =
c_{1} \ddot u_{0}(\tau) \left( R(\tau) - R(\bar\tau) \right) +
c_{1} \, {\rm e}^{-F(\tau)} /\dot u_{0}(\tau)$, so that
\begin{equation}
\dot w_{11}(0) = c_{1} \lim_{\tau\to 0} \left(
\ddot u_{0}(\tau)R(\tau) + \frac{{\rm e}^{-F(\tau)}}{\dot u_{0}(\tau)}
\right) - c_{1} \ddot u_{0}(0) R(\bar\tau) .
\label{eq:A13}
\end{equation}
Lemma \ref{lem:A.1} shows that the limit in (\ref{eq:A13})
is well defined, so that we obtain $\dot w_{11}(0)=0$ provided
\begin{equation}
R(\bar\tau) = \frac{1}{\ddot u_{0}(0)} \lim_{\tau\to 0}
\left( \ddot u_{0}(\tau) \, R(\tau) +
\frac{{\rm e}^{-F(\tau)}}{\dot u_{0}(\tau)} \right) .
\label{eq:A14}
\end{equation}
Since $R(\bar\tau)$ is finite, by (\ref{eq:A10}) also the
function $R(\tau)$ is strictly increasing in $\tau$
from $-\infty$ to $+\infty$. Therefore (\ref{eq:A14}) has
one and only one solution $\bar\tau$ in $(0,\pi\rho)$.\EP
%%%%%%%%%%%%%%%%%%%%%%%%%%%%%%%%%%%%%%%%%%%%%%%%%%%%%%%%%%%%%%%%%%%%%%%%%

%%%%%%%%%%%%%%%%%%%%%%%%%%%%%%%%%%%%%%%%%%%%%%%%%%%%%%%%%%%%%%%%%%%%%%%%%
%%%%%%%%%%%%%%%%%%%%%%%%%%%%%%%%%%%%%%%%%%%%%%%%%%%%%%%%%%%%%%%%%%%%%%%%%
\zerarcounters
\section{Non-vanishing of the constant $\boldsymbol{A}$}
\label{app:B}
%%%%%%%%%%%%%%%%%%%%%%%%%%%%%%%%%%%%%%%%%%%%%%%%%%%%%%%%%%%%%%%%%%%%%%%%%
%%%%%%%%%%%%%%%%%%%%%%%%%%%%%%%%%%%%%%%%%%%%%%%%%%%%%%%%%%%%%%%%%%%%%%%%%

Recall the definition (\ref{eq:8.3}) of $A$. We can write 
${\rm e}^{-f_{0}\tau} b(\tau)=w_{11}(\tau) - a(\tau) =
w_{11}(\tau) - \gamma \, \dot u_{0}(\tau)$, with
$\gamma = c_{2}/c$ and $\dot u_{0}(\tau) \, h(u_{0}(\tau)) +
2(\rho \Omega_{0})^{-2} k (u_{0}(\tau)) = - [ 2 \rho \Omega_{0}
\ddot u_{0}(\tau) + \dot u_{0}(\tau) \, h(u_{0}(\tau))]$
--- see (\ref{eq:6.7}), (\ref{eq:6.13}) and (\ref{eq:5.2}) ---,
so obtaining
\begin{equation}
A = - \frac{1}{2\pi \rho} \int_{0}^{2\pi \rho} {\rm d}\tau \,
{\rm e}^{F(\tau)} \left( w_{11}(\tau) - \gamma \, \dot u_{0}(\tau) \right)
\left[ 2 \rho \Omega_{0} \ddot u_{0}(\tau) +
\dot u_{0}(\tau) \, h(u_{0}(\tau)) \right] .
\label{eq:B4}
\end{equation}
%

%%%%%%%%%%%%%%%%%%%%%%%%%%%%%%%%%%%%%%%%%%%%%%%%%%%%%%%%%%%%%%%%%%%%%%%%%
\begin{lemma} \label{lem:B.1}
One has
$$ \int_{0}^{2\pi \rho} {\rm d}\tau \,
{\rm e}^{F(\tau)} \dot u_{0}(\tau)
\left[ 2 \rho \Omega_{0} \ddot u_{0}(\tau) +
\dot u_{0}(\tau) \, h(u_{0}(\tau)) \right] = 0 . $$
\end{lemma}
%%%%%%%%%%%%%%%%%%%%%%%%%%%%%%%%%%%%%%%%%%%%%%%%%%%%%%%%%%%%%%%%%%%%%%%%%

%%%%%%%%%%%%%%%%%%%%%%%%%%%%%%%%%%%%%%%%%%%%%%%%%%%%%%%%%%%%%%%%%%%%%%%%%
\prova By writing $\FF(\tau)={\rm e}^{F(\tau)}$,
one has $\dot\FF(\tau) = f(u_{0}(\tau))\,\FF(\tau) =
h(u_{0}(\tau))\,\FF(\tau)/ \rho\Omega_{0}$; cf. (\ref{eq:6.8}). Hence
$\FF \dot u_{0} \left[ 2 \rho \Omega_{0}\ddot u_{0} +
\dot u_{0} \, h(u_{0}) \right] =
\rho \Omega_{0} \left( \FF \frac{{\rm d}}{{\rm d}t}
\dot u_{0}^{2} + \dot \FF \dot u_{0}^{2} \right) =
\rho\Omega_{0} \frac{{\rm d}}{{\rm d}t} \left( \FF
\dot u_{0}^{2} \right)$, so that
\begin{eqnarray}
& & \null \hskip-1.5truecm
\int_{0}^{2\pi \rho} {\rm d}\tau \,
{\rm e}^{F(\tau)} \dot u_{0}(\tau) \left[ 2 \rho \Omega_{0} \ddot u_{0}(\tau) +
\dot u_{0}(\tau) \, h(u_{0}(\tau)) \right] =
\rho \Omega_{0} \int_{0}^{2\pi\rho} {\rm d}\tau \,
\frac{{\rm d}}{{\rm d}\tau} \left( \FF(\tau)\dot u_{0}^{2}(\tau) \right)
\nonumber \\
& & = \rho \Omega_{0} \left[
\FF(2\pi \rho)\dot u_{0}^{2}(2\pi \rho) -
\FF(0)\dot u_{0}^{2}(0) \right] =
\rho \Omega_{0} \left( \FF(2\pi \rho) -
\FF(0) \right) \dot u_{0}^{2}(0) = 0 ,
\label{eq:B6}
\end{eqnarray}
where we have used that $\dot u_{0}(\tau)$ is $2\pi\rho\,$-periodic
and $\dot u_{0}(0)=0$.\EP
%%%%%%%%%%%%%%%%%%%%%%%%%%%%%%%%%%%%%%%%%%%%%%%%%%%%%%%%%%%%%%%%%%%%%%%%% 

\vskip.2truecm

Because of Lemma \ref{lem:B.1}, (\ref{eq:B4}) becomes
\begin{equation}
A = - \frac{1}{2\pi \rho} \int_{0}^{2\pi \rho} {\rm d}\tau \,
{\rm e}^{F(\tau)} w_{11}(\tau)
\left[ 2 \rho \Omega_{0} \ddot u_{0}(\tau) +
\dot u_{0}(\tau) \, h(u_{0}(\tau)) \right] .
\label{eq:B7}
\end{equation}
%

%%%%%%%%%%%%%%%%%%%%%%%%%%%%%%%%%%%%%%%%%%%%%%%%%%%%%%%%%%%%%%%%%%%%%%%%%
\begin{lemma} \label{lem:B.2}
One has
$$ \int_{0}^{2\pi \rho} {\rm d}\tau \, {\rm e}^{F(\tau)}
w_{11}(\tau) \left[ 2 \rho \Omega_{0} \ddot u_{0}(\tau) +
\dot u_{0}(\tau) \, h(u_{0}(\tau)) \right] = -
2\pi \rho \, r_{1} \rho\Omega_{0} , $$
with $r_{1}$ defined in Remark \ref{rmk:2}.
\end{lemma}
%%%%%%%%%%%%%%%%%%%%%%%%%%%%%%%%%%%%%%%%%%%%%%%%%%%%%%%%%%%%%%%%%%%%%%%%%

%%%%%%%%%%%%%%%%%%%%%%%%%%%%%%%%%%%%%%%%%%%%%%%%%%%%%%%%%%%%%%%%%%%%%%%%%
\prova By writing once more $\FF(\tau)={\rm e}^{F(\tau)}$, we have
\begin{eqnarray}
& & \!\!\!\!\!\!\!\!\!\!\!\!
\FF \,w_{11} \left[ 2 \rho \Omega_{0} \ddot u_{0} +
\dot u_{0} \, h(u_{0}) \right] =
\rho \Omega_{0} w_{11} \left[
\FF \ddot u_{0} +
\left( \FF \ddot u_{0} + \dot \FF \dot u_{0} \right) \right] =
\rho \Omega_{0} \left[ w_{11} \FF \ddot u_{0} +
w_{11} \frac{{\rm d}}{{\rm d}\tau} \left( \FF\dot u_{0} \right)
\right] \nonumber \\
& & \!\!\!\!\!\!\!\!\!\!\!\!
= \rho \Omega_{0} \left[ w_{11} \FF \ddot u_{0} +
\frac{{\rm d}}{{\rm d}\tau} \left( \FF \dot u_{0} w_{11} \right)
- \FF \dot u_{0} \dot w_{11} \right] =
\rho \Omega_{0} \left[
\frac{{\rm d}}{{\rm d}\tau} \left( \FF \dot u_{0} w_{11} \right)
+ \FF \left( w_{11} \ddot u_{0} - \dot w_{11} \dot u_{0}
\right) \right] ,
\label{eq:B8}
\end{eqnarray}
where
\begin{equation}
\FF \left( w_{11} \ddot u_{0} - \dot w_{11} \dot u_{0}(\tau) \right) =
\frac{1}{c_{2}} \FF \left( w_{11} w_{22} - w_{21} w_{12}(\tau) \right) =
\frac{1}{c_{2}} \FF \det W = \frac{1}{c_{2}} \FF {\rm e}^{-F} =
\frac{1}{c_{2}} = r_{1} , \nonumber
\end{equation}
so that the integration of (\ref{eq:B8}) gives
\begin{eqnarray}
& & \null \hskip-1.5truecm \int_{0}^{2\pi \rho} {\rm d}\tau \,
{\rm e}^{F(\tau)} w_{11}(\tau) \left[ 2 \rho \Omega_{0} \ddot u_{0} +
\dot u_{0}(\tau) \, h(u_{0}(\tau)) \right] \nonumber \\
& & \null \hskip-1.truecm =
\rho \Omega_{0} \left[ \FF(2\pi \rho) \dot u_{0}(2\pi \rho)
w_{11}(2\pi \rho) - \FF(0) \dot u_{0}(0) w_{11}(0) +
2\pi \rho \, r_{1} \right] =
2\pi \rho \, r_{1} \rho\Omega_{0},
\label{eq:B10}
\end{eqnarray}
where once more we have used that $\dot u_{0}(2\pi \rho)=
\dot u_{0}(0)=0$.\EP
%%%%%%%%%%%%%%%%%%%%%%%%%%%%%%%%%%%%%%%%%%%%%%%%%%%%%%%%%%%%%%%%%%%%%%%%%

\vskip.2truecm

By using Lemma \ref{lem:B.2} in (\ref{eq:B7}) we obtain $A= - r_{1}
\rho\Omega_{0}$. Therefore $A \neq 0$ for any value $\rho\in\QQQ$.
Note that the time rescaling implies that $r_{1}$ is of the form
$r_{1}=(\rho\Omega_{0})^{-2}\bar r_{1}$, with $\bar r_{1}$
independent of $\rho$, so that $A=\bar A/\rho$, with
$\bar A=-\bar r_{1}/\Omega_{0}$ independent of $\rho$, consistently
with Remark \ref{rmk:4}.

%%%%%%%%%%%%%%%%%%%%%%%%%%%%%%%%%%%%%%%%%%%%%%%%%%%%%%%%%%%%%%%%%%%%%%%%%
%%%%%%%%%%%%%%%%%%%%%%%%%%%%%%%%%%%%%%%%%%%%%%%%%%%%%%%%%%%%%%%%%%%%%%%%%%
% References
%%%%%%%%%%%%%%%%%%%%%%%%%%%%%%%%%%%%%%%%%%%%%%%%%%%%%%%%%%%%%%%%%%%%%%%%%%
%%%%%%%%%%%%%%%%%%%%%%%%%%%%%%%%%%%%%%%%%%%%%%%%%%%%%%%%%%%%%%%%%%%%%%%%%

\end{document}